\documentclass[reqno]{amsart}

\usepackage{t1enc}
\usepackage[T2A]{fontenc}
\usepackage[cp1251]{inputenc}
\usepackage{amssymb, euscript, amsmath}
\usepackage{hyperref}

\begin{document}
\title[\hfilneg KARAZIN 2008 \hfil Thermoviscoelastic Berger plate]
{Asymptotic behaviour of thermoviscoelastic Berger plate}
\author[M. Potomkin\hfil KARAZIN 2008 \hfilneg]
{Mykhailo Potomkin}
\begin{abstract}
System of partial differential equations with a convolution terms and non-local nonlinearity describing oscillations of plate due to Berger's approach and with accounting for thermal regime in terms of Coleman-Gurtin and Gurtin-Pipkin law and fading memory of material is considered. The equation is transformed into a dynamical system in a suitable Hilbert space, which asymptotic behaviour is analysed. Existence of the compact global attractor in this dynamical system and some of its properties are proved in this article. Main tool in analysis of asymptotic behaviour is stabilizability inequality.     
\end{abstract}
\maketitle
\numberwithin{equation}{section}
\newtheorem{theorem}{Theorem}[section]
\newtheorem{proposition}[theorem]{Proposition}
\newtheorem{corollary}[theorem]{Corollary}
\newtheorem{lemma}[theorem]{Lemma}
\newtheorem{definition}[theorem]{Definition}
\newtheorem{remark}[theorem]{Remark}

\section {Introduction}
\normalsize 

Our main goal in this paper is to study long-time behaviour of the next system of integral-differential equations arising in plate theory  
$$
\left\{\begin{array}{ll}
\partial^{2}_{tt}u+k_1(0)\Delta^2u+\int\limits_{0}^{+\infty}k'_1(s)\Delta^2u(t-s)ds+
\left(\Gamma-\int_{\Omega}\left|\nabla u\right|^2d{\bf x}\right)\Delta u+ \nu\Delta v=p({\bf x})\\
\partial_t v-\omega\Delta v-(1-\omega)\int\limits_{0}^{+\infty}k_2(s)\Delta v(t-s)ds-\nu\Delta u_t=0,  
\end{array}\right.
$$
$$
{\bf x}=(x_1,x_2)\in\Omega\subset\mathbb R^2,\;t>0.
$$
with initial data 
$$
v(t,{\bf x})|_{t\leq 0}=v_0(-t,{\bf x}),\;u(t,{\bf x})|_{t\leq 0}=u_0(-t,{\bf x}). 
$$

Here we consider a thin plate of uniform thickness. When the plate is unloaded and is in null equilibrium its middle surface occupies a region $\Omega\subset\mathbb R^2$ of the plane $\left\{x_3=0\right\}$; $u(t,{\bf x})$ is a vertical component of displacement of corresponding point in middle surface. The presence of non-local term $\left(\Gamma-\int_{\Omega}\left|\nabla u\right|^2d{\bf x}\right)$, $\Gamma>0$ is explained by peculiarities in derivation of equation due to Berger's approach (see~\cite{Berger}). The first equation takes into account that material is viscous homogeneous and isotropic one, so convolution integral with the scalar kernel $k_1(s)$ appears (see~\cite{Renardy}). The function $v(t,{\bf x})$ is the temperature variation field and thus it satisfies one of the variant of heat equation. Here we consider heat equation according to Gurtin-Pipkin Law when $\omega=0$ (see~ \cite{Curt-Pipkin}) or Coleman-Gurtin Law when $\omega\in(0,1)$ (see \cite{ColeCurt}) instead of usual Furier Law when $\omega=1$, which has two main shortcomings. First, it is unable to account for the memory effects and, second, it predicts that a thermal disturbance at one point of the body is instantly felt everywhere in the body (for exact derivation of such heat equations for isotropic homogeneous material with memory see, e.g., \cite{GiorgiNaso,GiorgiPataA}). Parameter $\nu>0$ provides connection between deflection and temperature and depends on mechanical and thermal properties of the material (for more details see \cite{Lagnese}). 

Memory kernels $k_1(s)$ and $k_2(s)$ are supposed to be smooth decreasing convex functions and $k_2(s)$ vanishes at infinity, $k_1(\infty)>0$. 

In addition to equations and initial data we have to set boundary conditions following \cite{Lagnese}:
\begin{equation}\label{bnd1}
\begin{array}{l}
u=k_1(0)\Delta u+\int\limits_{0}^{+\infty}k'_{1}(s)\Delta u(t-s)ds=0,\;{\bf x}\in\partial\Omega,\;t\geq0,\\
v=0,\;{\bf x}\in\partial\Omega,\;t\in \mathbb R.
\end{array}
\end{equation}

These conditions are version of hinged boundary conditions simplified by the hypothesis that the action of boundary operator $B_1$ (for its definition and more details we refer to \cite{Lagnese,LasieTriggi}) is inessential and could be neglected. In this paper we provide condition (see Proposition \ref{propVolterra}) under which solutions of the considered problem satisfy more recognizable simplified hinged boundary conditions where the memory term is abscent (see, for example, \cite{ChuiBucci,Book,ChuiLasie,LinPata,TVE,LasieTriggi} and many others, where such boundary conditions were imposed for different models), namely,

\begin{equation} \label{bnd2}
u=\Delta u=0,\;{\bf x}\in\partial\Omega,\;t\geq0.
\end{equation}

To consider the model we will introduce new auxiliary variables which replace convolution integrals in original equation by some functional operator applied to one of the new added variable and allow us to apply the asymptotic theory of semigroups. Such approach originally being invented and applied in \cite{Dafermos} is widely-used in consideration of equations with memory (see \cite{SingPata,GiorgiPata,LinPata} and etc.). 

Linear versions of the model with memory in only thermal variable ($k'_1(s)\equiv 0)$ have been investigated in \cite{GiorgiPataA,LinPata}. Well-posedness, asymptotic stability, the presence and, in the same time, lack of exponential decay depending on conditions on thermal memory kernel were obtained in these works. Asymptotic stabilizability of a similar linear model but of the hyperbolic type when rotational forces are taken into account with clamped boundary conditions was considered in \cite{fabr}. Analogous work devoted to the linear thermoviscoelastic model has recently come out (see \cite{TVE}). Besides, questions of singular limit, i.e., asymptotic behaviour when kernels $k_i(s)$ tend to Dirac mass $\delta_0$ are considered in \cite{TVE}. Asymptotic behaviour (existence of compact global attactor) of homogeneous and isotropic viscoelastic solid described by semilinear hyperbolic equation was considered in \cite{ContiPata,ViscPata} without accounting for thermal regime. Models with memory are also investigated in \cite{gmpz,ppz}.

Isothermal Berger model of oscillations of plate without memory effects with the stress on its asymptotic behaviour was investigated in \cite{Book,Manuscript}. Up to our knowledge nonlinear model of the form considered in this paper with both viscoelastic and thermal memories was not studied before.  

Our main result is the proof of existence of compact global attractor of certain geometrical structure and of finite dimension. The proof is based on the method developed in \cite{ChuiLasie,Manuscript,LittleManu}, we refer also to \cite{ChuiBucci}. So-called stabilizability inequality (see Section 5) plays the crucial role in the proof. Such inequalities appeared in investigation of different kind of problems concerned with dissipative wave dynamics and becomes important tool in study of existence, smoothness and finite dimensionality of attractors (see \cite{Manuscript} and references therein). One should notice that these estimates are not consequences of some common abstract results and depend on peculiarities of the model under consideration in the essential degree. In slightly different form (from the one exploited in our paper) ideas of stabilizability inequality were developed in other works, e.g., \cite{ZME} and the method of $l$-trajectories developed in \cite{Malek1,Malek2}.
       
We conclude the Introduction with brief plan of the paper.

Section 2 is devoted to well-posedness result. In particular, after introducing all necessary settings the definition of a mild solution is given and then the question of its existence, uniqueness and dependence on initial data (Lipschitz property of the semiflow $S_t$) is considered. The Section includes the assertion of existence of classical solutions in the sense of semigroups. Besides, other questions like explicit representation formulas (Subsection 2.3), properties of the set of stationary points (Subsection 2.4), the existence of strict Lyapunov function (Subsection 2.5) complete the general (non-asymptotic) analysis of the semigroup. In addition results devoted to the differentiability of the semigroup and backward uniqueness (Subsection 2.6) end the Section 2. These results are needed for further asymptotic analysis.     

Section 3 includes main result, namely, the proof of existence of finite dimensional compact global attractor. It is divided by two parts. All necessary definitions and abstract results are given in Subsection 3.1. Subsection 3.2. includes the proof but the main part of it, namely, the proof of stabilizability estimate, is relegated to Section~ 5 because it is rather long and complicated and requires additional representational Lemma. Some properties of the attractor, in particular, obtained with the help of stabilizability inequality, are stated in Section~4. 

\section {Nonlinear Semigroup}
\subsection {Abstract form of the problem and main assumptions}
Let $\Omega$ be a bounded domain in $\mathbb R^2$ with smooth 
or rectangular boundary $\partial\Omega$, $\Delta$ denotes the Laplace operator.
We consider the following system of equations with linear memory 
\begin {equation}
\left\{
\begin{array}{l}
u_{tt}+k_1(0)\Delta^2u+\int\limits_{0}^{+\infty}k'_1(s)\Delta^2u(t-s)ds+\nu\Delta v=\\\;\;\;\;\;\;\;\;\;= p+ M\left(\int_{\Omega}\left|\nabla u\right|^2d{\bf x}\right)\Delta u,\\
v_t-\omega\Delta v-\int\limits_{0}^{+\infty}k_2(s)\Delta v(t-s)ds=\nu\Delta u_t,\\
u=k_1(0)\Delta u+\int\limits_{0}^{+\infty}k'_1(s)\Delta u(t-s)ds=0,\;\;{\bf x}\in\partial\Omega,\;t\geq0,\\
v=0,\;\;{\bf x}\in\partial\Omega,\;t\in\mathbb R\\      
u|_{t\leq0}=u_0(-t,{\bf x}),\;\;v|_{t\leq0}=v_0(-t,{\bf x}),\;\;{\bf x}\in\Omega.
\end{array}\right.
\label{1} \end{equation}

Now we intend to rewrite the system in abstract form, having replaced the Laplace operator defined on $H^2(\Omega)\cap H_{0}^{1}(\Omega)$ by an abstract self-adjoint positive operator $A$ which domain $D(A)$ is the subset of a Hilbert space $H$. 

Namely, we denote by $H$ a separable Hilbert space with inner product $(\cdot,\cdot)$ and corre\-spond\-ing norm $\left\|\cdot\right\|$. Let $A$ be a self-adjoint positive linear operator defined on a domain $D(A)\subset H$. Assume that there exists an eigenbasis $ \left\{e_k\right\}_{k=1}^{\infty}$ of operator $A$ such that
$$
(e_k,e_j)=\delta_{kj},\;\;Ae_k=\lambda_ke_k,\;\;k,j=1,2,...,
$$    
and 
$$
0<\lambda_1\leq\lambda_2\leq...,\;\;\lim\limits_{k\rightarrow\infty}\lambda_k=\infty,
$$
where $\lambda_k$ is corresponding eigenvalue of operator $A$.

We introduce the scale of Hilbert spaces $F_s$ in the next way
$$
F_s\equiv D(A^s)=\left\{v=\sum\limits_{k=1}^{\infty}c_ke_k\;:\;\sum\limits_{k=1}^{\infty}c_{k}^{2}\lambda_{k}^{2s}<\infty\right\},
$$ 
endowed with usual inner products:
$$
(v,w)_{s}=(A^sv,A^sw)=\sum\limits_{k=1}^{\infty}\lambda_{k}^{2s}(v,e_k)(w,e_k).
$$ 

As such $A:\;D(A)\subset H\rightarrow H$ we may take $A=-\Delta:\;H^2(\Omega)\cap H_{0}^{1}(\Omega)\subset L^2(\Omega)\rightarrow L^2(\Omega)$.

Next we replace kernels 
$$
\mu_1(s)=-k'_1(s),\;\;\mu_2(s)=-(1-\omega)k'_2(s).
$$
and we require 
\begin{eqnarray}
\mu_i(s)\in C^1(\mathbb R_+)\cap L^1(\mathbb R_+)\cap C[0,+\infty),\label{2}\\
\mu_i(s)\geq 0,\label{3}\\
\mu'_i(s)+\delta_i\mu_i(s)\leq0.\label{4}
\end{eqnarray}
where $\mathbb R_+=(0,+\infty)$.

Also we introduce weighted Hilbert spaces $L^{2}_{\mu_1}(\mathbb R_+;F_1)$ and $L^{2}_{\mu_2}(\mathbb R_+;F_{1/2})$ of measurable functions $\xi$ with values in $F_1$ or $F_{1/2}$ respectively such that
$$
\left\|\xi\right\|^{2}_{L^{2}_{\mu_1}(\mathbb R_+;F_{1})}\equiv
\int\limits_{0}^{+\infty}\mu_1(s)\left\|\xi(s)\right\|^{2}_{1}ds<\infty
$$
and 
$$
\left\|\xi\right\|^{2}_{L^{2}_{\mu_2}(\mathbb R_+;F_{1/2})}\equiv
\int\limits_{0}^{+\infty}\mu_2(s)\left\|\xi(s)\right\|^{2}_{1/2}ds<\infty.
$$

Following the ideas from \cite{Dafermos} we introduce additional variables, namely, the summed past history of $u$ and $v$, defined as 
$$
\overline{\eta}^t(s)=u(t)-u(t-s),\;\;\eta^{t}(s)=\int\limits_{0}^{s}v(t-y)dy,
$$  
they formally satisfy linear equations 
$$\begin{array}{l}
\frac{\partial}{\partial t}\overline{\eta}^{t}+\frac{\partial}{\partial s}\overline{\eta}^{t}= u_t(t)\\
\frac{\partial}{\partial t}{\eta}^{t}+\frac{\partial}{\partial s}{\eta}^{t}= v(t),
\end{array}
$$
and
$$
\overline{\eta}^t(0)=\eta^t(0)=0,
$$
whereas
$$
\overline{\eta}^0(s)=\overline{\eta}_0(s)\equiv u_0(0)-u_0(s),\;\;\;
\eta^{0}(s)=\eta_0(s)\equiv \int\limits_{0}^{s}v_0(y)dy. 
$$

The following Cartesian product of Hilbert spaces will play the role of a phase space for the considered model:
$$
\mathcal{H}=F_1\times F_0\times F_0 \times L^{2}_{\mu_1}(\mathbb R_+; F_1)\times L^{2}_{\mu_2}(\mathbb R_+;F_{1/2}) 
$$
with scalar product denoted as $\left\langle \cdot,\cdot \right\rangle$.

Let $\overline{T}$, $T$ be linear operators in $L^{2}_{\mu_1}(\mathbb R_+;F_1)$ and $L^{2}_{\mu_2}(\mathbb R_+;F_{1/2})$ respectively with domains 
$$
\begin{array}{c}
D(\overline{T})=\left\{\overline{\eta}\in L^{2}_{\mu_1}(\mathbb R_+;F_1)\left|\overline{\eta}_s\in L^{2}_{\mu_1}(\mathbb R_+;F_1),\;\overline{\eta}(0)=0\right.\right\}\\
D(T)=\left\{\eta\in L^{2}_{\mu_2}(\mathbb R_+;F_{1/2})\left|\eta_s\in L^{2}_{\mu_2}(\mathbb R_+;F_{1/2}),\;\eta(0)=0\right.\right\}
\end{array}
$$
defined by 
$$
\overline{T}\overline{\eta}=-\overline{\eta}_s, \;\;
T\eta=-\eta_s
$$
for all admissible $\overline{\eta}$ and $\eta$. Here $\eta_s$ denotes the distributional derivative with respect to the "memory" variable~$s$.   

These operators satisfy next inequalities 
$$\begin{array}{c}
\left(\overline{T}\overline{\eta},\overline{\eta}\right)_{L^{2}_{\mu_1}(\mathbb R_+;F_1)}\leq -\frac{\delta_1}{2}\left\|\overline{\eta}\right\|^{2}_{L^{2}_{\mu_1}(\mathbb R_+;F_1)},\;\;\forall \overline{\eta}\in D(\overline{T}),\\
\left( {T}{\eta},{\eta}\right)_{L^{2}_{\mu_2}(\mathbb R_+;F_{1/2})}\leq -\frac{\delta_2}{2}\left\|{\eta}\right\|^{2}_{L^{2}_{\mu_2}(\mathbb R_+;F_{1/2})},\;\;
\forall \eta\in D(T).
\end{array}
$$

We consider just first inequality. Its proof obtained with the help of integration by parts 
$$\begin{array}{l}
\int\limits_{0}^{+\infty}\mu_1(s)(-\frac{\partial}{\partial s}\overline{\eta}(s),\overline{\eta}(s))_1ds=
-\frac{1}{2}\int\limits_{0}^{+\infty}\mu_1(s)\frac{\partial}{\partial s}\left\|\overline{\eta}(s)\right\|_{1}^{2}ds=\\
=\frac{1}{2}\int\limits_{0}^{+\infty}\mu'_1(s)\left\|\overline{\eta}(s)\right\|_{1}^{2}ds\leq
-\frac{\delta_1}{2}\left\|\overline{\eta}\right\|^{2}_{L^{2}_{\mu_1}(\mathbb R_+;F_1)}
\end{array}
$$ 

Here we used requirements on the kernel. For more detailed argument see, e.g., \cite{SingPata, GiorgiPata} and references therein. 

For further investigations we are to impose conditions on function $M(\cdot)$, namely:
\begin{equation}
\left\{\begin{array}{l}
\mathcal{M}(z)\equiv\int\limits_{0}^{z}M(\xi)d\xi\geq-az-b,\;\;a\in(0,\lambda_1),\;b\in\mathbb R,\\
M(z)\in C^2(\mathbb R_+).
\end{array}\right.
\label{5}\end{equation}

In view of notation above \eqref{1} transforms into 
\begin {equation}
\left\{
\begin{array}{l}
u_{tt}+\beta A^2u+\int\limits_{0}^{+\infty}\mu_1(s)A^2\overline{\eta}^t(s)ds-\nu A v= p- M\left(\left\|A^{1/2}u\right\|^2\right)A u,\\
v_t+\omega A v+\int\limits_{0}^{+\infty}\mu_2(s)A{\eta}^t(s)ds+\nu A u_t=0,\\
\overline{\eta}^{t}_{t}=\overline{T}\overline{\eta}^{t}+u_t(t),\;\;\;\;
\eta^{t}_{t}=T\eta^{t}+v(t),\\
u|_{t=0}=u_0,\;\;u_t|_{t=0}=u_1,\;\;
v|_{t=0}=v_0,\;\;\overline{\eta}^t|_{t=0}=\overline{\eta}_0,\;\;
\eta^t|_{t=0}=\eta_0.
\end{array}\right.
\label{6}\end{equation}

The proof of existence and uniqueness is based on the theory of linear semigroups (see \cite{Pazy}). Therefore for the sake of convenience we represent linear part of equation \eqref{6} with the help of linear opertor $\mathcal{L}:D(L)\subset \mathcal{H}\rightarrow \mathcal{H}$ given by 
$$
{\mathcal{L}}U=\left(
\begin {array}{c}
w\\
-\beta A^2u-\int\limits_{0}^{ \infty }
\mu_1(s)A^2\overline{\eta}(s)ds+\nu A v\\-\omega A v-\int\limits_{0}^{\infty}\mu_2(s)A\eta(s)ds-\nu A w\\
 \overline{T}\overline{\eta}+w\\T{\eta}+v
\end{array}\right),
\;\;\;\;
U=\left(\begin{array}{c}
u\\w\\v\\\overline{\eta}\\{\eta}
\end{array}\right)\in\mathcal{H}.
$$         
and equipped with the domain:
$$
D(\mathcal{L})=\left\{U=\left(\begin{array}{c}
u\\w\\v\\\overline{\eta}\\{\eta}
\end{array}\right)\in\mathcal{H}\left|
\begin{array}{l}
\overline{\eta}\in D(\overline{T}),\;\eta\in D(T)\\
w\in F_1, \; v\in F_{1/2}\\
\beta A^2u + \int\limits_{0}^{+\infty} \mu_1(s)A^2\overline{\eta}(s)ds-\nu Av\in F_0\\
\omega Av+\int\limits_{0}^{+\infty}\mu_2(s)A\eta(s)ds\in F_0
\end{array}\right.\right\}
$$

In the next Section we prove that operator $\mathcal{L}$ is the infinitesimal operator of s.c. semigroup of contractions in space $\mathcal{H}$. 

Having made final notations for nonlinear term, namely,
$$f(U)=\left(
\begin{array}{c}
0\\-M\left(\left\|A^{1/2}u\right\|^2\right)A u+p\\0\\0\\0
\end{array}
\right),
$$
we rewrite nonlinear problem \eqref{6} as a first order problem of the form 
\begin{equation}
\left\{
\begin{array}{l}
\dot{U}(t)={\mathcal L}U(t)+f(U(t))\\
U(0)=U_0 \in \mathcal{H}
\end{array}
\right.
\label{7}\end{equation}

We recall that according to \cite{Pazy} $U(t)$ is a {\it mild solution} of \eqref{7} if $U(t)$ satisfies the following equality
$$
U(t)=e^{t\mathcal{L}}U_0+\int\limits_{0}^{t}e^{(t-\tau)\mathcal{L}}f(U(\tau))d\tau,
$$
where $e^{t\mathcal{L}}$ is the linear semigroup on $\mathcal{H}$ which infinitesimal operator is $\mathcal{L}$. $U(t)$ is called a {\it classical solution} on interval $[0,T)$ if it is continuously differentiable, its values lie in $D(\mathcal{L})$ and it satisfies \eqref{7}.   

\subsection {Generation of Semigroup.}
\normalsize 

In this Section we prove well-posedness result formulated in the Theorem below. The proof consists of several steps. First, the problem with only linear part exploiting the notion of infinitesimal operator is considered. Then according to corresponding Theorems from \cite{Pazy} existence and uniqueness result is obtained. In addition, there are assertions devoted to continuous dependence on initial data and the existence of classical solutions in the formulation of the Theorem. Together they yield that solutions of the problem \eqref{7} generate continuous semigroup of non-linear operators according to definition from \cite{Book}.
             
\begin{theorem} \label{2.2.1} 
\it Let assumptions \eqref{2},\eqref{3},\eqref{4} and \eqref{5} hold true. Assume also that $p\in H$. Then for all $U_0 \in \mathcal{H}$ and $T>0$ there exists a unique mild solution $U(t)\in C(0,T;\mathcal{H})$. 
 
Besides, if $U_1,U_2\in \mathcal{H}$ and $\left\|U_i\right\|_{\mathcal{H}}\leq R$ then there exists a positive constant $C_{R,T}$ such as 
\begin{equation}
\left\|S_tU_1-S_tU_2\right\|_{\mathcal{H}}\leq C_{R,T}\left\|U_1-U_2\right\|_{\mathcal{H}},\;\;t\in\left[0,T\right].
\label{8}\end{equation}

And if $U_0\in D(\mathcal{L})$ then the corresponding mild solution $U(t)$ is a classical solution.  
\end{theorem}
\rm 
{\it Proof.}

STEP I. In order to prove that $\mathcal{L}$ defined in the previous Subsection is the infinitesimal of s.c. semigroup of contractions we use Lumer-Phillips Theorem (see \cite{Pazy}), thus, we need to show $\mathcal{L}$ to be maximal and dissipative one. For similar arguments see \cite{SingPata,GiorgiPataA,LinPata,TVE,Lagnese,Renardy}.

The property of being a dissipative one, i.e.
$$
<\mathcal{L}U,U>_{\mathcal{H}}\leq0\;\;\forall U\in D(\mathcal{L}),
$$ is obvious if one redefine the norm of $\mathcal{H}$ and equipped scalar product into equivalent one, via
$$\left\|U\right\|^{2}_{\mathcal{H}}=\beta\left\|Au_0\right\|^2+
\left\|w\right\|^2+\left\|v\right\|^2+\left\|\overline{\eta}\right\|^{2}_{L^{2}_{\mu_1}(\mathbb R_+;F_1)}+\left\|\eta\right\|^{2}_{L^{2}_{\mu_2}(\mathbb R_+;F_{1/2})}  $$

The operator $\mathcal{L}$ is the maximal one provided that the mapping $I-\mathcal{L}:D(\mathcal{L})\rightarrow\mathcal{H}$ is onto. Let $U^{\ast}=(u^{\ast};w^{\ast};v^{\ast};\overline{\eta}^{\ast};\eta^{\ast})\in\mathcal{H}$, and consider the equation
$$
(I-\mathcal{L})U=U^{\ast}
$$        
which, written in components, reads 
\begin{eqnarray}
&&u-w=u^{\ast}\in F_1 \label{9}\\
&&w+\beta A^2u+\int\limits_{0}^{+\infty}\mu_1(s)A^2\overline{\eta}(s)ds-\nu Av=w^{\ast}\in F_0\label{10}\\
&&v+\omega Av+\int\limits_{0}^{+\infty}\mu_2(s)A\eta(s)ds+\nu Aw=v^{\ast}\in F_0\label{11}\\
&&\overline{\eta}+\overline{\eta}_s-w=\overline{\eta}^{\ast}\in L^{2}_{\mu_1}(\mathbb R_+;F_1)\label{12}\\
&&\eta+\eta_s-v=\eta^{\ast}\in L^{2}_{\mu_2}(\mathbb R_+;F_{1/2}) \label{13}
\end{eqnarray}

Integrations of two latter equalities immediately implies that
\begin{eqnarray}
&&\overline{\eta}(s)=w(1-e^{-s})+\int\limits_{0}^{s}e^{y-s}\overline{\eta}^{\ast}(y)dy \label{14}\\
&&\eta(s)=v(1-e^{-s})+\int\limits_{0}^{s}e^{y-s}\eta^{\ast}(y)dy. \label{15}
\end{eqnarray} 

Sabsituting \eqref{14} and \eqref{15} into \eqref{10} and \eqref{11} respectively, accounting for  
$$
\int\limits_{0}^{+\infty}\mu_1(s)A^2\int\limits_{0}^{s}e^{y-s}\overline{\eta}^{\ast}(y)dy\in F_{-1},\;\;
\int\limits_{0}^{+\infty}\mu_2(s)A\int\limits_{0}^{s}e^{y-s}\eta^{\ast}(y)dy\in F_{-1/2},
$$
we reduce original system \eqref{9}-\eqref{13} to the system of three equations 
$$
\begin{array}{l}
u-w=u^{\ast}\in F_1\\
w+\beta A^2u+c_1A^2w-\nu Av=w^{\ast\ast}\in F_{-1}\\
v+\omega Av+c_2Av+\nu Aw=v^{\ast\ast}\in F_{-1/2}
\end{array}
$$
where elements $w^{\ast\ast}$ and $v^{\ast\ast}$ are supposed to be given and  
$$
c_i=\int\limits_{0}^{+\infty}\mu_i(s)(1-e^{-s})ds,\;\;i=1,2.
$$

Or it could be rewritten in terms of only $w,\;v$ as follows 
\begin{eqnarray}
&&w+c_{\beta}A^2w-\nu Av=w^{\ast\ast\ast}\in F_{-1} \label{16}\\
&&v+c_{\omega}Av+\nu Aw=v^{\ast\ast\ast}\in F_{-1/2}, \label{17} 
\end{eqnarray}
where equations are obtained by substitution the relation $u=w+u^{\ast}$ into the latter system, elements $w^{\ast\ast\ast}$, $v^{\ast\ast\ast}$ are also supposed to be given, $c_\beta$ and $c_\omega$ are positive constants.  

To solve the elliptic problem \eqref{16}-\eqref{17} we apply Lax-Millgram Theorem with settings like in \cite{Lions}. Namely,

$$V=F_1\times F_{1/2},\;\mathrm{H}=F_0\times F_0,\;V^{\ast}=F_{-1}\times F_{-1/2}$$
$$
\begin{array}{lcr}
a((w,v);(\tilde{w},\tilde{v}))&=&(w,\tilde{w})+c_{\beta}(Aw,A\tilde{w})-\nu(Av,\tilde{w})+\\
&&+(v,\tilde{v})+c_{\omega}(Av,\tilde{v})+\nu(Aw,\tilde{v}).\end{array}$$

$V^{\ast}$ being the dual of $V$ with respect to $\mathrm{H}$ and the bilinear form $a((w,v);(\tilde{w},\tilde{v}))$ being coercetive, Lax-Millgram Theorem is applicable and implies the existence of $w\in F_{1}$ and $v\in F_{1/2}$ that satisfy \eqref{16}-\eqref{17}. The element $U=\left(u;w;v;\overline{\eta};\eta\right)$, where $u=u^{\ast}+w$ and "memory"\; components - $\overline{\eta}$ and $\eta$ - are obtained by \eqref{14} and \eqref{15}, satisfies the system of equations \eqref{9}-\eqref{13} and so - on account for the form of these equalities - obviously belongs to $D(\mathcal{L})$. Thus $\mathcal{L}$ is the maximal operator and due to Lumer-Phillips Theorem generates s.c. semigroup of contractions. 

STEP II. The existence of local solutions is the consequence of \cite[Theorem 6.1.4]{Pazy}. More precisely, $\forall U_0\in \mathcal{H}$ $\exists t_{max}\leq\infty$ and there exists a unique function $U(t)\in C\left([0,t_{max});\mathcal{H}\right)$ such as $U(t)$ is the mild solution of \eqref{7} on each closed interval $[0,T]$ where $T<t_{max}$. Besides, if $t_{max}<\infty$ then
\begin{equation}
\lim\limits_{t\uparrow t_{max}}\left\|U(t)\right\|_{\mathcal{H}}=\infty.
\label{18}\end{equation}
 
Naturally, application of this Theorem is allowed because each of its conditions is satisfied. Namely, linear part of the problem $-$ the operator $\mathcal{L}$ $-$ is generator of s.c. semigroup and nonlinearity $-$ function $f(U)$ $-$ is locally Lipschitz one. 
The statement that any mild solution could be extended to arbitrary closed interval of the form $[0,T]$ is equivalent to the equality $t_{max}=\infty$. 

Consider any mild solution $U(t)$ with initial data $U_0$. Assume that $t_{max}<\infty$ and $0<T<t_{max}$. Hence, \eqref{18} takes place. Next we apply \cite[Theorem 4.2.7]{Pazy}. According to this Theorem there exist sequences $\left\{f_n(t)\right\}_{n=1}^{\infty}\subset C^1\left([0,T];\mathcal{H}\right)$ and $\left\{U_{0n}\right\}_{n=1}^{\infty}\in D(\mathcal{L})$ such as 
$$
\begin{array}{l}
f_n(t)\rightarrow f(U(t))\;\mathrm{in}\; L^1(0,T;\mathcal{H})\\
U_{0n}\rightarrow U_0\;\mathrm{in}\;\mathcal{H}
\end{array}
$$   

Besides, there exists a sequence $\left\{U_n(t)\right\}_{n=1}^{\infty}$ of functions that satisfy next Coushy problem
$$
\left\{\begin{array}{l}
\frac{dU_n}{dt}=\mathcal{L}U_n(t)+f_n(t), \;t\in [0,T]\\
U_n(0)=U_{0n}.
\end{array}\right.
$$

Then for $\forall T'<T$ the sequence of $U_n(t)$ converges to $U(t)$ uniformly for all $t\in[0,T']$. 

Moreover, the following inequality holds true 
$$
\left\|U_n(t)\right\|^{2}_{\mathcal{H}}-\left\|U_{0n}\right\|^{2}_{\mathcal{H}}\leq 2 \int\limits_{0}^{t}\left\langle f_n(s),U_n(s)\right\rangle_{\mathcal{H}}ds.
$$ 

Passing to the limit $n\rightarrow+\infty$ we obtain 
\begin{equation}
\left\|U(t)\right\|^{2}_{\mathcal{H}}-\left\|U_{0}\right\|^{2}_{\mathcal{H}}\leq 2 \int\limits_{0}^{t}\left\langle f(U(s)),U(s)\right\rangle_{\mathcal{H}}ds.
\label{19}\end{equation}

Using the same procedure considering the equality 
$$
\mathcal{M}\left(\left\|A^{1/2}u_n(t)\right\|^2\right)-
\mathcal{M}\left(\left\|A^{1/2}u_{0n}\right\|^2\right)=
\int\limits_{0}^{t}M\left(\left\|A^{1/2}u_n(s)\right\|^2\right)\left(Au_n(s),w_n(s)\right)ds,
$$
where $u_n$ and $w_n$ are corresponding components of $U_n$, we obtain  
\begin{equation}
\mathcal{M}\left(\left\|A^{1/2}u(t)\right\|^2\right)-
\mathcal{M}\left(\left\|A^{1/2}u_0\right\|^2\right)=
\int\limits_{0}^{t}M\left(\left\|A^{1/2}u(s)\right\|^2\right)\left(Au(s),w(s)\right)ds.
\label{20}\end{equation}

Next we consider the sum of \eqref{19} and \eqref{20}. Before this we set 
$$
E(t)\equiv\frac{1}{2}\left\|U(t)\right\|^{2}_{\mathcal{H}}+\mathcal{M}\left(\left\|A^{1/2}u\right\|^2\right)
$$ 

Using conditions on $\mathcal{M}(\cdot)$ and $p$, we obtain the next chain of inequalities 
$$
\alpha_1\left\|U(t)\right\|^{2}_{\mathcal{H}}-C_1\leq E(t)\leq E(0)+\int\limits_{0}^{t}(p,u(s))ds\leq C_2+\alpha_2\int\limits_{0}^{t}\left\|U(\tau)\right\|^{2}_{\mathcal{H}}d\tau.  
$$
Here and below all new constants are positive.

Then
$$
\left\|U(t)\right\|^{2}_{\mathcal{H}}\leq C\left(\left\|U_0\right\|^{2}_{\mathcal{H}}\right)+\alpha\int\limits_{0}^{t}\left\|U(\tau)\right\|^{2}_{\mathcal{H}}d\tau.
$$

Application of Gronwall Lemma is left:
\begin{equation}
\left\|U(t)\right\|^{2}_{\mathcal{H}}\leq C\left(\left\|U_0\right\|^{2}_{\mathcal{H}}\right)e^{C_3t}.
\label{21}\end{equation} 

That obviously contradicts to \eqref{18}. Thus we have proved that every mild solution could be extended on a closed interval of arbitrary length.   

STEP III. We continue the proof considering the question of continuous dependence of the solution on initial data. 

Consider $\forall T>0,\;\forall t\in(0,T)$ and two mild solutions $U_1(t)$ and $U_2(t)$ with initial data $U_{10}$ and $U_{20}$ respectively, then 
$$
\left\|U_1(t)-U_2(t)\right\|_{\mathcal{H}}\leq\left\|e^{t\mathcal{L}}(U_{10}-U_{20})\right\|_{\mathcal{H}}+\int\limits_{0}^{t}\left\|e^{(t-\tau)\mathcal{L}}\left(f(U_1(\tau))-f(U_2(\tau))\right)\right\|_{\mathcal{H}}d\tau 
$$  

Using that $\left\|e^{t\mathcal{L}}\right\|_{[\mathcal{H},\mathcal{H}]}\leq 1$, estimate \eqref{21} and locally Lipschitz property of $f$ with corresponding constant $L(R)$ (i.e., $f$ is the Lipschitz function in the closed ball $\left\{\left\|U\right\|^2\leq R\right\}$ with constant $L(R)$, here it is reasonable to set $$R\equiv C(\max\left\{\left\|U_{10}\right\|^{2}_{\mathcal{H}},\left\|U_{20}\right\|^{2}_{\mathcal{H}}\right\})e^{C_3T}$$ where all constants are taken from \eqref{21}) and again Gronwall Lemma we finally obtain 
$$
\left\|U_1(t)-U_2(t)\right\|_{\mathcal{H}}\leq e^{C_T}\left\|U_{10}-U_{20}\right\|_\mathcal{H} 
$$  
where $C_T$ is a positive constant that depends on initial data. 

STEP IV. The statement about classical solutions follows directly from \cite[Theorem 6.1.5]{Pazy}.
  
The proof is complete.

Now we may set  
$
S_tU_0\equiv U(t)
$, then $(\mathcal{H},S_t)$ is the dynamical system on $\mathcal{H}$ that is generated by mild solutions of \eqref{7} (for exact definition of a dynamical system see \cite{BabinVishik,Book,Temam}).

We continue with observation that is of interest in its own rights and not used in asymptotic analysis. In what folows below in this Subsection we will impose conditions on initial data from domain of operator $\mathcal{L}$ under which the corresponding classical solution (having returned to original problem with settings $H=L^2(\Omega),$ $A=-\Delta,$ $D(A)=H^2(\Omega)\cap H^{1}_{0}(\Omega)$) satisfies boundary conditions \eqref{bnd2}. 

Necessity of additional conditions to satisfy \eqref{bnd2} is illustrated by the next example. 

Consider $U=(u;w;v;\overline{\eta};\eta)\in D(\mathcal{L})$ given as follows
$$
\begin{array}{rl}
u=\sum\limits_{k\geq1}\frac{1}{k\lambda_k}e_k, & \overline{\eta}(s)=-s\left(\frac{\beta}{\overline{\kappa}_1}u-\frac{\nu}{\overline{\kappa}_1}\sum\limits_{k\geq1}\frac{1}{k\lambda_{k}^{3/2}}e_k\right),\\
v=\sum\limits_{k\geq1}\frac{1}{k\lambda_{k}^{1/2}}e_k,& \eta(s)=-\frac{\omega}{\overline{\kappa}_2}sv,
\end{array}
$$    
where $\overline{\kappa}_i=\int_{0}^{\infty}s\mu_i(s)ds$ and we recall that $e_k$ and $\lambda_k$ is corresponding eigenvector and eigenvalue of operator $A$ respectively. The component $w\in F_1$ is arbitrary. 

In this case, in particular, 
$$\begin{array}{c}
\int\limits_{0}^{+\infty}\mu_2(s)A\eta(s)ds \notin F_0, \;
\int\limits_{0}^{+\infty}\mu_1(s)A^2\overline{\eta}(s)ds \notin F_0.
\end{array}
$$  

And $\beta A^2u-\nu Av\notin F_0.$ Hence, $Au\notin F_{1/2}$. We recall that in terms of original problem \eqref{1} $F_{1/2}=H^{1}_{0}(\Omega).$ Therefore conditions \eqref{bnd2} does not hold. 

The main difficulty is in the fact that we may conclude that the sum 
\begin{equation}\label{diez}
\beta A^2u+\int\limits_{0}^{+\infty}\mu_1(s)A^2\overline{\eta}(s)ds 
\end{equation}   
lies in the space $F_{-1/2}$ but we can't say the same separately for each part of this sum.  

Nevertheless, it turned out that if we impose additional conditions on initial data we will manage to separate two parts in \eqref{diez}. Namely, next Proposition takes place. 

\begin{proposition} \label{propVolterra}
Let $U_0\in D(\mathcal{L})$ and, moreover, 
$$
u_0\in L^{\infty}(0,+\infty;F_{3/2}),
$$
where $u_0(t)\equiv u_0-\overline{\eta}_0(t)$ for all $t\geq0$. 
Then the corresponding classical solution satisfies 
\begin{equation}\label{VolterraCond}
u\in C([0,+\infty);F_{3/2}).
\end{equation}

And, hence, if, moreover, $H=L^2(\Omega),$ $A=-\Delta,$ $D(A)=H^2(\Omega)\cap H^{1}_{0}(\Omega)$ then 
$$
\Delta u(t,{\bf x})=0,\;\;{\bf x}\in\partial\Omega,\;t\geq0.
$$
\end{proposition}

{\it Proof.}
First we note that $\overline{\eta}^t(s)=u(t)-u(t-s)$ where $u(-t)=u_0(t),\;t\geq0$. For general case of a mild solution this formula will be proved in Subsection 2.3. but one can see that we just returned to introduction of the memory variable in Subsection~2.1. 

Next, from equations \eqref{6} and the formula above for $\overline{\eta}$ we obtain 
\begin{equation}\label{Volterra1}
u(t)-\int\limits_{-\infty}^{t}\frac{\mu_1(t-y)}{\kappa_1+\beta}u(y)dy=h(t),\;t\geq0. 
\end{equation}  
where $\kappa_1\equiv\int\limits_{0}^{+\infty}\mu_1(s)ds$ and $h(t)\in C([0,\infty);F_{3/2})$. To obtain injection for $h(t)$ that satisfies $$(\kappa_1+\beta)A^2 h(t)=-u_{tt}+\nu Av-M\left(\left\|A^{1/2}u\right\|^2\right)Au+p$$ one should use Theorem \ref{dif} for continuity of derivatives $u_{tt}$, $v_t$ and $\eta^{t}_{t}$ and manner of the proof of the estimate \eqref{forv} in the Corollary \ref{smooth2} for continuity of $v(t)$ with values in $F_{1/2}$). 

Equation \eqref{Volterra1} may be rewritten
\begin{equation}\label{Volterra2}
u(t)-\int\limits_{0}^{t}\frac{\mu_1(t-y)}{\kappa_1+\beta}u(y)dy=F(t),\;\;t\geq 0,
\end{equation}   
where 
$$
F(t)=h(t)+\int\limits_{-\infty}^{0}\frac{\mu_1(t-y)}{\kappa_1+\beta}u_0(-y)dy.
$$

Note that $F(t)$ belongs to $C([0,+\infty);F_{3/2})$. We will solve \eqref{Volterra2} by standard iteration method on interval $[0,T]$ where $T>0$ is arbitrary. Namely, we set $w_0=0$, 
$$
w_n(t)=F(t)+\int\limits_{0}^{t}\frac{\mu_1(t-y)}{\kappa_1+\beta}w_{n-1}(y)dy,\;\;n=1,2,...,
$$       
and we observe that 
$$
\sup\limits_{[0,T]}\left\|w_{n+1}(t)-w_n(t)\right\|_{3/2}\leq q \cdot \sup\limits_{[0,T]}\left\|w_{n}(t)-w_{n-1}(t)\right\|_{3/2} \leq q^n \cdot \sup\limits_{[0,T]}\left\|F(t)\right\|_{3/2}
$$
where 
$$
q=\frac{\kappa_1}{\kappa_1+\beta}<1.
$$ 

Thus $\left\{w_n(t)\right\}$ is a Cauchy sequence in $C([0,T];F_{3/2})$, and it converges to $u(t)\in C([0,T];F_{3/2})$. For last conclusion we need to say that solution of \eqref{Volterra2} is unique since the operator 
$$
\int\limits_{0}^{t}\frac{\mu_1(t-y)}{\kappa_1+\beta}\bullet dy\;:C([0,T];F_{3/2})\rightarrow C([0,T];F_{3/2})
$$ 
is an operator of contractions.

The proof is complete. 

\begin{remark} 
Though the manner of solvation of Volterra equation \eqref{Volterra2} is standard it should be noted that similar equations in study of viscous models were considered in \cite{Dafermos,Lagnese}. 
\end{remark}

\subsection{Explicit representation formula.}
In the sequel we need typical for equations with infinite memory explicit representation formulas (similar to considered in \cite{ContiPata,SingPata,fabr,GiorgiPataA,ViscPata,LinPata}).
  
\begin{proposition}\label{331}
Let $U(t)=(u(t);w(t);v(t);\overline{\eta}^t;\eta^t)$ be a mild solution of \eqref{7} with initial data $U_0=(u_0;w_0;v_0;\overline{\eta}_0;\eta_0)$. Then
\begin{equation}
\overline{\eta}^t(s)=
\left\{\begin {array}{ll}
u(t)-u(t-s),&t>s>0\\
\overline{\eta}_0(s-t)+u(t)-u(0),&t\leq s
\end{array}
\right.
\label{24}\end{equation} 
\end{proposition}

\begin{proposition}\label{332}
 Let $U(t)=(u(t);w(t);v(t);\overline{\eta}^t;\eta^t)$ be a mild solution of \eqref{7} with initial data $U_0=(u_0;w_0;v_0;\overline{\eta}_0;\eta_0)$. Then
\begin{equation}
\eta^t(s)=
\left\{\begin {array}{ll}
\int\limits_{0}^{s}v(t-y)dy,&t>s>0\\
\eta_0(s-t)+\int\limits_{0}^{t}v(t-y)dy,&t\leq s
\end{array}
\right.
\label{25}\end{equation}
\end{proposition}

{\it Proof.}

\rm We restrict ourselves to the case of Proposition \ref{331}. Other Proposition is proved in the same manner.   
First we note that each mild solution of \eqref{7} could be approximated by classical solutions of the problem. More precisely, for all $U_0\in\mathcal{H}$ we can choose sequence $\left\{U_{0n}\;:\;U_{0n}\in D(\mathcal{L})\right\}$ such as $U_{0n}\rightarrow U_0$ in $\mathcal{H}$ (such choise is possible since $D(\mathcal{L})$ is dense in $\mathcal{H}$) and due to Theorem \ref{2.2.1} for arbitrary $T>0$: 
$$
\left.\begin{array}{l}
\exists U_n(t)\;\mathrm{a\;classical\;solution\;of}\;\mathrm{\eqref{7}}\\
\exists U(t)\;\mathrm{a\;mild\;solution\;of}\;\mathrm{\eqref{7}}
\end{array}
\right|\;U_n(t)\rightarrow U(t)\;\mathrm{uniformly\;on}\;\mathrm{[0,T]} 
$$ 

Here we present the derivation of explicit representation formulas (of course, reader can just verify formulas substituting them into corresponding equations in \eqref{6}).
Now we derive explicit representation formula for the first "memory"  component of the classical solution $U_n(t)=(u_n(t);u_{t,n}(t);v_n(t);\overline{\eta}^{t}_{n};\eta^{t}_{n})$ . Consider the third equation of system \eqref{6}:
$$
\frac{\partial}{\partial t} \overline{\eta}^{t}_{n}(s)=-
\frac{\partial}{\partial s} \overline{\eta}^{t}_{n}(s)+u_{t,n}(t)
$$

Then after the substitution $y=t-s$ we obtain 
$$
\frac{\partial}{\partial t} \overline{\eta}^{t}_{n}(t-y)=
\frac{\partial}{\partial y} \overline{\eta}^{t}_{n}(t-y)+u_{t,n}(t)
$$

And in account for 
$
\frac{d}{dt}\overline{\eta}^{t}_{n}(t-y)=
\frac{\partial}{\partial t}\overline{\eta}^{t}_{n}(t-y)-
\frac{\partial}{\partial y}\overline{\eta}^{t}_{n}(t-y)
$
we obtain
$$
\frac{d}{dt}\overline{\eta}^{t}_{n}(t-y)=u_{t,n}(t)
$$

To reach the final equality the process of integration is left:
$$
\begin{array}{lll}
\mathrm{let}\;t>s,\;\mathrm{integration}\;\left.\int\limits_{y}^{t}\cdot\right|&&
\overline{\eta}^{t}_{n}(t-y)-\overline{\eta}^{y}_{n}(0)=u_n(t)-u_n(y)\\
&\mathrm{or}&\overline{\eta}^{t}_{n}(s)=u_n(t)-u_n(t-s)\\
\mathrm{let}\;t\leq s\;\mathrm{integration}\;\left.\int\limits_{0}^{t}\cdot\right|&&
\overline{\eta}^{t}_{n}(t-y)-\overline{\eta}^{0}_{n}(-y)=u_n(t)-u_n(0)\\
&\mathrm{or}&\overline{\eta}^{t}_{n}(s)=\overline{\eta}_{0,n}(s-t)+u_n(t)-u_n(0). 
\end{array}
$$

We used above that $U_{0n}\in D(\mathcal{L})$ (and it implies that $\overline{\eta}^{y}_{n}(0)=0$) and initial condition (namely, $\overline{\eta}^{0}_{n}(-y)=\overline{\eta}_{0,n}(-y)$).

Our next step is typical. To obtain necessary equalities for $U(t)$ we pass to limit $n\rightarrow\infty$. Before this we denote 
$$
\psi^{t}_{n}(s)=\left\{
\begin{array}{ll}
u_n(t)-u_n(t-s),&t>s>0\\
\overline{\eta}_{0,n}(s-t)+u_n(t)-u_n(0),&t\leq s
\end{array}\right.$$and$$
\psi^{t}(s)=\left\{
\begin{array}{ll}
u(t)-u(t-s),&t>s>0\\
\overline{\eta}_{0}(s-t)+u(t)-u(0),&t\leq s
\end{array}\right..
$$   

We have already known that $\psi^{t}_{n}(s)=\overline{\eta}^{t}_{n}(s)$. We need $\psi^{t}(s)=\overline{\eta}^{t}(s)$. 

Since $\overline{\eta}^{t}_{n}(s)\rightarrow\overline{\eta}^t(s)$ in $L^{2}_{\mu_1}(\mathbb R_+;F_1)$ uniformly on $t\in\mathrm{[0,T]}$, it is sufficient to show that $\psi^{t}_{n}(s)\rightarrow\psi^{t}(s)$ in $L^{2}_{\mu_1}(\mathbb R_+;F_1)$ for all $t\in\mathrm{[0,T]}$. 

Indeed, consider any $t\in\mathrm{[0,T]}$:
$$
\left\|\psi_{n}^{t}-\psi^{t}\right\|_{L^{2}_{\mu_1}(\mathbb R_+;F_1)}^{2}=
\int\limits_{0}^{+\infty}\mu_1(s)\left\|\psi_{n}^{t}(s)-\psi^{t}(s)\right\|^{2}_{1}ds=
$$$$=\int\limits_{0}^{t}\mu_1(s)\left\|\left(u_n(t)-u(t)\right)-\left(u_n(t-s)-u(t-s)\right)\right\|^{2}_{1}ds+
$$    
$$
+\int\limits_{t}^{+\infty}\mu_1(s)\left\|\left(u_n(t)-u(t)\right)-\left(u_n(0)-u(0)\right)+\left(\overline{\eta}_{0,n}(s-t)-\overline{\eta}_0(s-t)\right)\right\|_1ds\rightarrow0.
$$   

Thus we may conclude $\overline{\eta}^{t}(s)=\psi^{t}(s)$ and this completes the proof.


\subsection{The set of stationary points.}

In this Subsection we analyse the set of stationary points of the problem \eqref{7} 
$$
\left\{
\begin{array}{l}
\dot{U}(t)={\mathcal L}U(t)+f(U(t))\\
U(0)=U_0 \in \mathcal{H},
\end{array}
\right.
$$
which could be defined as follows
$$
\mathcal{N}=\left\{U\in X\;:\;S_tU=U \;\forall\; t\geq0\right\}. 
$$ 

We note that stationary point $U_0\in\mathcal{H}$ is the mild solution of \eqref{7} $U(t)\equiv U_0$ and, as a consequence, it satisfies the following integral equation
$$
U_0=e^{t\mathcal{L}}U_0+\int\limits_{0}^{t}e^{(t-\tau)\mathcal{L}}f(U_0)d\tau.
$$    

This yields that for any $t>0$ 
$$
-\left(\frac{e^{t\mathcal{L}}-I}{t}\right)U_0=\frac{1}{t}\int\limits_{0}^{t}e^{(t-\tau)\mathcal{L}}f(U_0)d\tau=\frac{1}{t}\int\limits_{0}^{t}e^{\tau\mathcal{L}}f(U_0)d\tau.
$$

Right-hand side converges to $f(U_0)$ as $t\downarrow0$ (see \cite[Theorem 1.2.4.(a)]{Pazy}). 
Therefore, by the definition of infinitesimal generator $U_0\in D(\mathcal{L})$ and
\begin{equation}
\mathcal{L}U_0+f(U_0)=0. 
\label{22}\end{equation} 

Thus we have next assertion:

\begin{proposition} The set $\mathcal{N}$ of stationary points could be written as follows:

\begin{equation}
\mathcal{N}=\left\{V=\left(u;0;0;0;0\right)\;:\;\beta A^2 u +M\left(\left\|A^{1/2}u\right\|^2\right)Au=p\right\}
\label{23}\end{equation} 
\rm
\end{proposition}

Properties of the set \eqref{23} when $\beta>0$ was investigated in \cite{Book}. In particular, boundedness of $\mathcal{N}$ was proved and conditions which implies finiteness of $\mathcal{N}$ were obtained. In general, results concerning the set $\mathcal{N}$ could be stated as follows (see \cite[Chapter 4]{Book})  
  
\begin{theorem} \label{ChuN} Let $\mathcal{J}[u]\equiv\beta A^2 u +M\left(\left\|A^{1/2}u\right\|^2\right)Au$ and $\mathcal{J}'[u]$ is its Freshet derivative for $u\in F_0$. We introduce the set 
$$\mathcal{R}\equiv\left\{h\in F_0\;:\;\exists\left[\mathcal{J}'[u]\right]^{-1}\;\;\mathrm{for}\;\mathrm{all}\;u\in\mathcal{J}^{-1}[h]\right\}$$

Then 
\begin{list}{}{}
\item[{(i)}] for any bounded $B\subset F_0$ preimage $\mathcal{J}^{-1}(B)$ is bounded (in particular, $\mathcal{N}$ is bounded in $\mathcal{H}$)
\item[{(ii)}] the set $\mathcal{R}$ is open, dense in $F_0$ and if $p\in\mathcal{R}$ then $\mathcal{N}$ is a finite set.  
\end{list}
\end{theorem} 
\rm

It should be noted that if a property of a dynamical system holds for the parameters from an open and dense set in the corresponding space, then it its frequently said that this property is a {\it generic property}. Generic properties are frequently encountered and stay stable during the small perturbations of the properties of a system (see \cite[Chapter 2]{Book}).
 
For illustration we consider the case when $M(z)=z-\Gamma$ and $p=0$ that corresponds to genuine (non-abstract) homogeneous Berger's equation. This case is described by the next statement that is easy to verify.

\begin{proposition}
Each stationary point has the form of $U=\left(u;0;0;0;0\right)$ where 
$$
u=c_k e_k, \;\;k=0,\pm 1,\pm 2,...,\pm N_0,
$$ 
$e_k$ $-$ eigenbasis vector of the operator $A$, $N_0$ is the maximal integer such that $\Gamma>\beta\lambda_{N_0}$ and 
$$
\begin{array}{l}
c_0=0,\\
c_{\pm k}=\pm\sqrt{\frac{\Gamma-\beta\lambda_k}{\lambda_k}}\;\;k=\overline{1,N_0}.
\end{array}
$$ 
\end{proposition}

\rm
Obviously this Proposition implies that the number of stationary points in considering case (and, we recall that $\beta>0$) is finite. 

\subsection {Strict Lyapunov function.}
 
It turned out that the semigroup $(\mathcal{H},S_t)$ which we consider in this work is gradient (see definition below). This circumstance alows to simplify asymptotic analysis due to well-known results (see Subsection 3.1). 

\begin{definition} The dynamical system $(X,S_t)$ is said to be {\rm gradient} if it possesses a strict Lyapunov function, i.e. there exists a continuous functional $\Phi(U)$ defined on $X$ such that (i) the function $t\rightarrow \Phi(S_tU)$ is nonincreasing for any $U\in X$, and (ii) the equation $\Phi(S_tU)=\Phi(U)$ for all $t>0$ implies that $S_tU=U$ for all $t>0$, i.e., $U$ is a stationary point of $(X,S_t)$.
\end{definition}   

Corresponding functional has the following form:
$$
\Phi(U)=\frac{1}{2}\left\|U\right\|_{\mathcal{H}}^{2}+\mathcal{M}\left(\left\|A^{1/2}u\right\|^2\right)-(p,u)
$$

Now we notice that each classical solution satisfies the energy relation 
\begin{eqnarray}
&\Phi(U(t))-\Phi(U(\tau))=-\omega\int\limits_{\tau}^{t}\left\|A^{1/2}v\right\|^2dy+\nonumber\\&+\int\limits_{\tau}^{t}\left( \overline{T}\overline{\eta}^{y},\overline{\eta}^{y}\right)_{L^{2}_{\mu_1}(\mathbb R_+;F_1)} dy+
\int\limits_{\tau}^{t}\left(T\eta^{y},\eta^{y}\right)_{L^{2}_{\mu_2}(\mathbb R_+;F_{1/2})} dy. 
\label{26}\end{eqnarray}

Therefore, for any mild solution we have the estimate 
\begin{equation}
\Phi(U(t))-\Phi(U(\tau))\leq-
\int\limits_{\tau}^{t}\left\| \overline{\eta}^{y}\right\|^{2}_{L^{2}_{\mu_1}(\mathbb R_+;F_1)} dy-
\int\limits_{\tau}^{t}\left\|\eta^{y}\right\|^{2}_{L^{2}_{\mu_2}(\mathbb R_+;F_{1/2})} dy.
\label{26m}\end{equation} 

The (energy) relation \eqref{26m} with Propositions \ref{331} and \ref{332} gives us the following result:

\begin{theorem}
 Let the functional $\Phi(U)\;:\;\mathcal{H}\longmapsto \mathbb R$ is given by 
$$
\Phi(U)=\frac{1}{2}\left\|U\right\|_{\mathcal{H}}^{2}+\mathcal{M}\left(\left\|A^{1/2}u\right\|^2\right)-(p,u)
$$  

Then 

I. The system $(\mathcal{H},S_t)$ is gradient with $\Phi$ as a Lyapunov function, i.e.

\begin{list}{}{}
\item[{(i)}] the function $t\mapsto\Phi(S_tU_0)$ is nonincreasing for any $U_0\in\mathcal{H}$;
\item[{(ii)}] the equation $\Phi(S_tU_0)=\Phi(U_0)$ for all $t>0$ and for some $U_0\in\mathcal{H}$ implies that $U_0$ is a stationary point. 
\end{list}        

II. The functional $\Phi(U)$ is bounded from above on any bounded subset of $\mathcal{H}$ and the set $\Phi_R=\left\{U\;:\;\Phi(U)\leq R\right\}$ is bounded for every $R$. 

Thus, $\Phi(U)$ is a appropriate strict Lyapunov function. 
\end{theorem}
\rm 

The statement I.(i) is proved with the help of relation \eqref{26m}, I.(ii) needs explicit representation formulas (Propositions \ref{331} and \ref{332}) besides \eqref{26m}. Statements in II hold true thanks to conditions imposed on function $\mathcal{M}$ and their proof requires manipulations the same as in proof of global existence (see Theorem \ref{2.2.1}, step II) so it is omitted here.

\subsection {Some other useful properties.}

Here we collect some more statements about the considered semigroup. We note that the statement devoted to Frechet differentiability of $S_t$ is similar to \cite[Proposition 2.3]{ChuiLasie} and backward uniquiness result for thermoelastic plates was obtained also in \cite{ChuiLasie}, but the case with memory variables is much simplier, what is noted in \cite{SingPata}. 

Consider the system that could be obtained after formal differentiation with the respect to $t$ of \eqref{7}
\begin{equation}\label{differ}
\left\{
\begin{array}{l}
\dot{W}=\mathcal{L}W+f'(U(t))W,\\
W(0)=W_0.
\end{array}
\right.
\end{equation}     

Here for $U(t)=(u(t);u_t(t);v(t);\overline{\eta}^t;\eta^t)$ and $W(t)=(w(t);w_t(t);\xi(t);\tilde{\overline{\eta}}^t;\tilde{\eta}^t)$
$$
f'(U(t))W=\left(\begin{array}{c}
0\\-M'(\left\|A^{1/2}u\right\|^2)(Au,w)Au-M(\left\|A^{1/2}u\right\|^2)Aw\\0\\0\\0\end{array}\right)^{\mathrm{T}}
$$

Using the standard method presented in this Section well-posedness result for \eqref{differ} is proved on the phase space $\mathcal{H}$ and moreover (compare with \eqref{21})
\begin{equation} \label{bnd}
\left\|W(t)\right\|_{\mathcal{H}}\leq e^{a_{R,T}}\left\|W_0\right\|_{\mathcal{H}},\;\;t\in[0,T]  
\end{equation}
provided $\left\|U(t)\right\|_{\mathcal{H}}\leq R$ for all $t\in [0,T]$.

Denote also 
$$
\begin{array}{l}
B(u)=p-M(\left\|A^{1/2}u\right\|^2)Au,\\
B'(u)w=-M'(\left\|A^{1/2}u\right\|^2)(Au,w)Au-M(\left\|A^{1/2}u\right\|^2)Aw.
\end{array}
$$

\begin{theorem} \label{dif}
The mapping $U\rightarrow S_tU$ is Frechet differentiable on $\mathcal{H}$ for every $t\geq 0$. Moreover, the Frechet derivative $D[S_tU_0]:\mathcal{H}\rightarrow\mathcal{H}$ is a mapping of the form 
\begin{equation}
D[S_tU_0]W_0=W(t)=(w(t);w_t(t);\xi(t);\tilde{\overline{\eta}}^t;\tilde{\eta}^t),\;\;
W_0=(w_0;w_1;\xi_0;\tilde{\overline{\eta}}_0;\tilde{\eta}_0), 
\end{equation}
where $(w(t);w_t(t);\xi(t);\tilde{\overline{\eta}}^t;\tilde{\eta}^t) \in C([0,\infty);\mathcal{H})$ is a unique solution to the problem \eqref{differ}.
\end{theorem}

{\it Proof.}
Consider $U_0,W_0\in\mathcal{H}$, $t\geq 0$ and the function 
$$
Y(t)=S_t[U_0+W_0]-S_t[U_0]-W(t).
$$

We need to show that
\begin{equation} \label{td}
 \left\|Y(t)\right\|_{\mathcal{H}}=\overline{O}(\left\|W_0\right\|_{\mathcal{H}}).
\end{equation}

Note that $Y(t)$ solves 
$$
\left\{
\begin{array}{l}
\dot{Y}=\mathcal{L}Y+\mathcal{F}(t),\\
Y(0)=0.
\end{array}
\right.
$$ 
where second component of $\mathcal{F}(t)$ (we denote it as $F(t)$, other components are equal to zero) is equal to 
$$
F(t)=B(u^*(t))-B(u(t))-B'(u(t))w(t),
$$ 
where $u^*(t)$, $u(t)$ and $w(t)$ are first components of $S_t[U_0+W_0]$, $S_t[U_0]$ and $W(t)$ respectively. The first component of $Y(t)$ will be denoted by $z(t)$. 

Next representation holds 
$$
F(t)=\mathrm{I}_1+\mathrm{I}_2,
$$
where 
$$
\begin{array}{l} 
\mathrm{I}_1=\int\limits_{0}^{1}\left[B'(u_{\lambda}(t))-B'(u(t))\right]w(t)d\lambda=\\ 
=-\int\limits_{0}^{1}\left\{
\left[M'(\left\|A^{1/2}u_{\lambda}\right\|^2)-M'(\left\|A^{1/2}u\right\|^2)\right](Au_{\lambda},w)Au_{\lambda}+\right.\\\\+M'(\left\|A^{1/2}u\right\|^2)(A(u_{\lambda}-u),w)Au_{\lambda}+
M'(\left\|A^{1/2}u\right\|^2)(Au,w)A(u_{\lambda}-u)+\\\\+\left.\left[M'(\left\|A^{1/2}u_{\lambda}\right\|^2)-M'(\left\|A^{1/2}u\right\|^2)\right]Aw\right\}d\lambda
\end{array}
$$
and
$$
\mathrm{I}_2=\int\limits_{0}^{1}B'(u_{\lambda}(t))z(t)d\lambda. 
$$
where $u_{\lambda}=u+\lambda(u^*-u)$. Henceforth we assume that all functions - $S_t[U_0~+~W_0]$, $S_t[U_0]$ and $W(t)$ - are bounded on $[0,T]$ with respect to the norm of $\mathcal{H}$ with number~$R$.

Using \eqref{21} and \eqref{bnd} we obtain 
$$
\left\|\mathrm{I}_1\right\|\leq C_R\left\|u^*(t)-u(t)\right\|_1\left\|w(t)\right\|_1\leq C_R\left\|W_0\right\|^{2}_{\mathcal{H}}.
$$

From energetical equation of the problem for $Y(t)$ we obtain 
$$
\left\|Y(t)\right\|^{2}_{\mathcal{H}}-\left\|Y(0)\right\|^{2}_{\mathcal{H}}\leq\int\limits_{0}^{t}(F(\tau),z_t)d\tau\leq C_R\left\|W_0\right\|^{4}_{\mathcal{H}}+ C_R\int\limits_{0}^{t}\left\|Y(\tau)\right\|^{2}_{\mathcal{H}}d\tau  
$$  

The final conclusion follows from Gronwall Lemma 
$$
\left\|Y(t)\right\|_{\mathcal{H}}\leq C_R\left\|W_0\right\|^{2}_{\mathcal{H}}.
$$

The proof is complete. 

Other additional result states injectivity of $S_t$ and of its Frechet derivative $D[S_tU_0]$ for all $t>0$ and $U_0\in \mathcal{H}$. Due to finite memory we can easily obtain the result which will be needed in Subsection 4.3.  

\begin{proposition}
Next statements hold:
\begin{list}{}{}
\item[$\bullet$] Let $$U_i(t)=(u^{i}(t);u^{i}_{t}(t);v^i(t);\overline{\eta}^{i,t};\eta^{i,t}),\;\;i=1,2$$  be two solutions of \eqref{7}. 

If $U^1(T)=U^2(T)$ for some $T>0$, then $U^1(t)=U^2(t)$ for every $t~\in~[0,T]$. 
\item[$\bullet$] Let $u(t)\in C([0,T];F_1)$ and $W(t)=$$(w(t);w_t(t);\xi(t);\tilde{\overline{\eta}}^t;\tilde{\eta}^t)$ be a solution to the linear (non-autonomous) equation \eqref{differ}.

If $W(T)=0$, then $W(t)=0$ for every $t\in[0,T]$.     
\end{list}\end{proposition} 

{\it Proof.} 
For pair of solutions both \eqref{7} and \eqref{differ} explicit representation formulas formulated in Propositions \ref{331} and \ref{332} hold. Therefore further proof is general for both problems. 

We have 
$$
\overline{\eta}^{1,T}(s)=\overline{\eta}^{2,T}(s)\;\;\forall s\geq0.
$$  

Then 
$$
u^1(T)-u^1(T-s)=u^2(T)-u^2(T-s) \;\;\forall s\in[0,T]
$$  

In view that $u^1(T)=u^2(T)$ it means
$$
u^1(t)=u^2(t)\;\;\forall t\in[0,T].
$$

Because of same arguments $v^1(t)=v^2(t)$ for all $t\in[0,T]$. 

The fact that memory variables coincide in initial moment is left to verify. It follows from the next representation 
$$
\overline{\eta}^{i,t}(s)=\overline{\eta}^{i,T}(s+T-t)-u^i(T)+u^i(t),\;\;t\in[0,T],\;i=1,2. 
$$
and the similar for $\eta^t(s)$.  

The proof is complete. 

\section {Main result: existence of finite dimensional attractor.}
\normalsize

\subsection{Preliminaries and formulation of main result.}

\rm \normalsize
Now we recall some definitions and statements (following mostly \cite{BabinVishik,Book,Temam}) that will be needed in the sequel. All formulations are made for abstract dynamical system $(X,S_t)$ where $X\;-$ is a metric space and $S_t$ is a semigroup of operators in $X$.  

\begin{definition} $\mathcal{A}\subset X$  is called an attractor if (i) $\mathcal{A}$ is closed bounded strictly invariant set ($S_t\mathcal{A}=\mathcal{A}\;\forall t\geq0$) and (ii) $\mathcal{A}$ possesses the uniform attraction property, i.e. for any bounded set $B\subset X$ the following equality holds true 
$$
\lim\limits_{t\rightarrow +\infty}\sup\limits_{U\in B}\mathrm{dist}_X \left(S_tU,\mathcal{A}\right)=0.
$$\rm
\end{definition}
    \rm 
    
\begin{definition} The dynamical system $(X,S_t)$ is said to be { \rm asymptotically smooth} if for any positively invariant bounded set $D\subset X$ there exists a compact $K$ in the closure $\overline{D}$ of $D$ such that 
$$
\lim\limits_{t\rightarrow +\infty}\sup\limits_{U\in D}\mathrm{dist}_X \left(S_tU,K\right)=0.
$$
\end{definition}\rm

To prove the existence of compact global attractor we rely on the following well-known assertion (see \cite{Manuscript,Hale}), that is useful in our case because it requires dynamical system to be gradient what has already been proved in the previous Section. Other advantage of this approach is abscence of necessity to obtain dissipativity first.  

\begin{theorem} Assume that $(X,S_t)$ is a gradient dynamical system which, moreover, is asymptotically smooth. Assume that Lyapunov function $\Phi(U)$ associated with the system is bounded from above on any bounded subset of  $X$  and the set $\Phi_R=\left\{U:\Phi(U)\leq R\right\}$ is bounded for every $R$. If the set $\mathcal{N}$ of stationary points of $(X,S_t)$ is bounded, then $(X,S_t)$ possesses a compact global attractor. 
\label{31}\end{theorem} 
\rm 
It turns out that in our case of a gradient system thanks to well-known statements (see \cite{BabinVishik,Book,ChuiLasie,Temam}) it is possible to describe geometrical structure of the attractor. 

\begin{definition} We define the {\rm unstable manifold} $\mathcal{M}^{u}\left(\mathcal{N}\right)$ emanating from the set $\mathcal{N}$ as a set of all $U\in X$ such that there exists a full trajectory $\gamma=\left\{U(t)\;:\;t\in\mathbb R\right\}$ with the properties 
$$
U(0)=U\;\;\mathrm{and}\;\;\lim\limits_{t\rightarrow-\infty}\mathrm{dist}_X(U(t),\mathcal{N})=0.
$$ 
\end{definition}\rm

The following assertion describes a long-time behaviour in terms of unstable manifold when the power of the set $\mathcal{N}$ (finite or infinite) is not specified. 

\begin{theorem} \label{geo}
Assume that the gradient system $(X,S_t)$ possesses a compact global attractor $\mathcal{A}$. Then $\mathcal{A}=\mathcal{M}^u\left(\mathcal{N}\right)$ and, moreover,  
\begin{list}{}{}
\item[{(i)}] the global attractor $\mathcal{A}$ consists of full trajectories $\gamma=\left\{U(t)\;:\;\mathbb R\right\}$ such that 
$$
\lim\limits_{t\rightarrow-\infty}\mathrm{dist}_X(U(t),\mathcal{N})=0\;\mathrm{and}\;
\lim\limits_{t\rightarrow+\infty}\mathrm{dist}_X(U(t),\mathcal{N})=0.
$$  
\item[{(ii)}] for any $U\in X$ we have 
$$
\lim\limits_{t\rightarrow+\infty}\mathrm{dist}_X(S_tU,\mathcal{N})=0. 
$$
\end{list}
\end{theorem}

Thus if all conditions of the Theorem above are satisfied then any trajectory stabilizes to the set $\mathcal{N}$ of stationary points. Assumption that $\mathcal{N}=\left\{e_1,...,e_n\right\}$ - is a finite set allows us to describe asymptotic behaviour more precise. Namely, next direct consequence of previous Theorem holds true:

\begin{corollary} \label{finit}
Assume that the gradient dynamical system $(X,S_t)$ possesses a compact global attractor $\mathcal{A}$ and $\mathcal{N}=\left\{e_i|\;i=\overline{1,n},\;e_i\in X\right\}$ is a finite set. Then $\mathcal{A}=\cup_{i=1}^{n}\mathcal{M}^u(e_i)$ and 
\begin {list}{}{}
\item[{(i)}] the global attractor $\mathcal{A}$ consists of full trajectories $\gamma=\left\{U(t):t\in\mathbb R\right\}$ connecting pairs of stationary points, i.e. any $U\in \mathcal{A}$ belongs to some full trajectory $\gamma$ and for any $\gamma\subset \mathcal{A}$ there exists a pair $\left\{e,e^{\ast}\right\}\subset\mathcal{N}$ such that 
$$
U(t)\rightarrow e\; \mathrm{as}\; t\rightarrow-\infty\;\mathrm{and}\;U(t)\rightarrow e^{\ast} \;\mathrm{as}\;t\rightarrow+\infty;   
$$  
\item[{(ii)}] for any $V\in X$ there exists a stationary point $e$ such that $S_tV\rightarrow e$ as $t\rightarrow+\infty.$
\end{list}
\end{corollary}
\rm

Therefore to obtain an existence of compact global attractor of the certain geometrical structure we have to investigate questions that concern with the set of stationary points, existence of a strict Lyapunov function and asymptotically smoothness of considered semigroup.
First two questions have already been considered in the previous Section. So we need to prove just asymptotically smoothness of the dynamical system $(\mathcal{H},S_t)$.  

An important characteristic of a global attractor is its dimension. We use here generalisation of notion "dimensionality". Namely, 

\begin{definition} The \rm fractal dimension \it $\mathrm{dim}^{X}_{f}M$ of a compact set $M$ in a complete metric space $X$ is defined by 
$$
\mathrm{dim}_{f}^{X}M=\limsup\limits_{\varepsilon\rightarrow 0}\frac{\ln N(M,\varepsilon)}{\ln(1/\varepsilon)},
$$ 
where $N(M,\varepsilon)$ is the minimal number of closed sets in $X$ of the diameter $2\varepsilon$ which cover the set $M$. \end{definition}
\rm

The proof of finite dimensionality is based on the next abstract result which is generalization of the Ladyzhenskaya's Theorem on the dimension of the invariant sets. To see  examples of application of this Theorem we refer to, e.g., \cite{ChuiBucci,ChuiLasie,Manuscript}. 

\begin{theorem} \label{Ladfinite}
Let $X$ be a Banach space and $M$ be a bounded closed set in $X$. Assume that there exists a mapping $V\;:\;M\mapsto X$ such that $M\subseteq VM$ and also 
\begin {list}{}{}
\item[(i)] V is Lipschitz on M, i.e., there exists $L>0$ such that 
$$
\left\|Vv_1-Vv_2\right\|\leq L\left\|v_1-v_2\right\|,\;\;v_1,v_2\in M; 
$$ 
\item[(ii)] there exist compact seminorms $n_1(x)$ and $n_2(x)$ on $X$ such that 
$$
\left\|Vv_1-Vv_2\right\|\leq\eta\left\|v_1-v_2\right\|+K\left[n_1(v_1-v_2)+n_2(Vv_1-Vv_2)\right]
$$ 
for any $v_1,v_2\in M$, where $0<\eta<1$ and $K>0$ are constants (a seminorm $n(x)$ on $X$ is said to be compact if for any bounded set $B\subset X$ there exists a sequence $\left\{x_n\right\}\subset B$ such that $n(x_n-x_m)\rightarrow 0$ as $m,n\rightarrow\infty$).
\end {list} 

Then $M$ is a compact set in a $X$ of a finite fractal dimension. Moreover, we have the estimate 
$$
\mathrm{dim}_{f}^{X}M\leq\left[\ln\frac{2}{1+\eta}\right]^{-1}\cdot\ln m_0\left(\frac{4K(1+L^2)^{1/2}}{1-\eta}\right),
$$    
where $m_0(R)$ is the maximal number of pairs $(x_i,y_i)$ in $X\times X$ possessing the properties 
$$
\left\|x_i\right\|^2+\left\|y_i\right\|^2\leq R^2,\;\;n_1(x_i-x_j)+n_2(y_i-y_j)>1,\;\;i\neq j.
$$\rm
\end{theorem}
  
Now we may formulate the main result of this section:  

\begin{theorem} \label{3A}
Assume that conditions \eqref{2},\eqref{3},\eqref{4},\eqref{5} and $p\in H$ hold. Then the dynamical system $(\mathcal{H},S_t)$ possesses a compact global atractor of the form $\mathcal{A}=\mathcal{M}^{u}(\mathcal{N})$ of finite fractal dimension.  
\end{theorem}
\rm 

\subsection{Proof of Theorem \ref{3A}.}
 
The following criterium (see \cite{CeronLopes,Manuscript}) leads to desired property (asymptotical smoothness):

\begin{theorem} \label{341}
Let $(X,S_t)$ be a dynamical system on a complete metric space $X$ endowed with a metric $d$. Assume that for any bounded positively invariant set $B$ in $X$ there exist numbers $T>0$ and $0<q<1$, and a pseudometric $\rho^{T}_{B}$ on $C(0,T;X)$ such that
\begin {list}{}{}
\item [(i)] the pseudometric $\rho^{T}_{B}$ is precompact (with respect to $X$) in the following sense: any sequence $\left\{x_n\right\}\subset B$ has a subsequence $\left\{x_{n_k}\right\}$ such that the sequence $\left\{y_k\right\}\subset C(0,T;X)$ of elements $y_k(\tau)=S_{\tau}x_{n_k}$ is Couchy with respect to $\rho^{T}_{B}$;   
\item [(ii)] the following inequality holds 
$$
d(S_Ty_1,S_Ty_2)\leq q\cdot d(y_1,y_2)+\rho^{T}_{B}(\left\{S_{\tau}y_1\right\},\left\{S_{\tau}y_2\right\}),\;\;
$$
for every $y_1,y_2\in B$, where we denote by $\left\{S_{\tau}y_i\right\}$ the element in the space $C(0,T;X)$ given by function $y_i(\tau)=S_{\tau}y_i$. 
\end{list}       
Then $(X,S_t)$ is an asymptotically smooth dynamical system. 
\end{theorem}

\rm  
Reader is refered to \cite[Chapter 2]{Manuscript} for details and other relative statements.

To apply the criterium above we obtain so-called "stabilizability inequality"  stated in the next Theorem. This Theorem will be proved in Section 5.   

\begin{theorem}  \label{342}
Assume $M(z)\in C^2(\mathbb R_+)$.
Let $(u^1;v^1;\overline{\eta}^1;\eta^1)$ and $(u^2;v^2;\overline{\eta}^2;\eta^2)$ be two solutions of the problem 
\eqref{7} with initial data $U^i=(u^{i}_{0};u^{i}_{1};v^{i}_{0};\overline{\eta}^{i}_{0};\eta^{i}_{0}),\;i=1,2$. Assume that 
$$\begin{array}{ll}
\left\|Au^{i}(t)\right\|^2+\left\|u_{t}^{i}(t)\right\|^2+\left\|v^{i}(t)\right\|^2+\left\|\overline{\eta}^{i,t}\right\|_{L^{2}_{\mu_1}(\mathbb R_+;F_1)}^{2}+
\left\|\eta^{i,t}\right\|_{L^{2}_{\mu_2}(\mathbb R_+;F_{1/2})}^{2}\leq R^2
& 
\end{array}$$
for all $t\geq0$. Let 
$$
Z(t)\equiv\left(u^1(t)-u^2(t);u_{t}^{1}(t)-u^{2}_{t}(t);v^1(t)-v^2(t);\overline{\eta}^{1,t}-\overline{\eta}^{2,t};\eta^{1,t}-\eta^{2,t}\right)
$$  
and 
$$
z(t)\equiv u_1(t)-u_2(t).
$$

Then there exist positive constants $C_R$ and $\gamma$ such that 

\begin{equation}
\left|Z(t)\right|^2\leq C_R\left|Z(0)\right|^2 e^{-\gamma t}+C_R \sup\limits_{0\leq\tau\leq t}
\left\|z(\tau)\right\|^2.
\label{27}\end{equation}
\end{theorem}
\rm

Now to apply both Theorem \ref{342} and \ref{341} we set 
$$
\begin{array}{l}
t\equiv T,\\
\rho_{B}^{T}(\left\{S_{\tau}y_1\right\},\left\{S_{\tau}y_2\right\})\equiv C_R\max\limits_{\tau \in [0,T]}\left\|u^1(\tau)-u^2(\tau)\right\|,\\
q\equiv C_Re^{-\gamma T}<1.
\end{array}
$$

Since $C(0,T;F_1)\cap C^1(0,T;H)$ compactly imbedded in $C(0,T;H)$ (see for example \cite{Simon}), pseudometric $\rho_{B}^{T}$ is precompact. Thus by Theorem \ref{341} $(\mathcal{H},S_t)$ is an asymptotically smooth dynamical system.

Therefore it follows from Theorems \ref{31} and \ref{geo} the compact global attractor $\mathcal{A}$ exists and possesses the structure of unstable manifold $\mathcal{A}=\mathcal{M}^{u}\left(\mathcal{N}\right)$. 

But Theorem \ref{3A} also asserts finite dimensionality of $\mathcal{A}$. For the complete proof of this assertion with the same stabilizability inequality immanented to the equation under consideration (but with other phase space that does not essentially change the proof) we refer to \cite{ChuiLasie,Manuscript} or discusssion in \cite{ChuiBucci}. 

To prove finiteness of the fractal dimension, we appeal to a generalization of the Ladyzhenskaya's Theorem on the dimension of the invariant sets (see Theorem \ref{Ladfinite}). This result applicable, in view of the local Lipschitz continuity of the semi-flow $S_t$ (see \eqref {8}) and of the stabilizability estimate.

Following the method described in \cite{ChuiLasie}, let us introduce the extended space $\mathcal{H}_{T}=\mathcal{H}\times W_1(0,T)$ (with an appropriate $T>0$). Here 
$$
W_1(0,T)=\left\{z(t)\;:\;\left|z\right|^{2}_{W_1(0,T)}\equiv\int\limits_{0}^{T}(\left\|Az(t)\right\|^2+\left\|z_t(t)\right\|^2)dt<\infty\right\}.
$$          

Next, we consider in $\mathcal{H}_T$ the set 
$$
\mathcal{A}_T:=\left\{U\equiv(u(0);u_t(0);v(0);\overline{\eta}^0;\eta^0;u(t),t\in[0,T]):(u(0);u_t(0);v(0);\overline{\eta}^0;\eta^0)\in \mathcal{A} \right\}, 
$$
where $$(u(t);u_t(t);v(t);\overline{\eta}^t;\eta^t)$$ is the solution to \eqref{7} with initial data $(u(0);u_t(0);v(0);\overline{\eta}^0;\eta^0)$, and define operator $V\;:\;\mathcal{A}_T\mapsto\mathcal{H}_T$ by the formula 
$$
V:(u(0);u_t(0);v(0);\overline{\eta}^0;\eta^0)\mapsto(u(T);u_t(T);v(T);\overline{\eta}^T;\eta^T;u(T+t)). 
$$

Then, by using pretty much the same arguments as in \cite{ChuiLasie,Manuscript}, we see that assumptions of Theorem \ref{Ladfinite} are satisfied. 

Thus proof of Theorem \ref{3A} is complete.

\section{Other properties of asymptotic behaviour.}

\subsection {Smoothness of the attractor.}
Often it's possible to prove that an attractor is the bounded set with respect to more strong topology (see for example \cite{ChuiBucci,ChuiLasie,ContiPata}). In order to obtain similar property for our case we use stabilizability estimate along with full invariance property of $\mathcal{A}$ like in \cite{ChuiLasie}. Besides, peculiarities of considered problem requires additional steps in order to obtain sufficiently explicit estimates.

First let us denote as $R>0$ such positive constant that 
\begin{equation} \label{87}
\left\|U_0\right\|_{\mathcal{H}}\leq R,\;\;\forall U_0\in \mathcal{A}.
\end{equation}

Our main goal in this Subsection is to prove step by step that there exists a positive constant $C_R$ such that for any trajectory $U(t)=(u(t);u_t(t);v(t);\overline{\eta}^t;\eta^t)$ lying in the attractor we have
\begin{eqnarray}
&\left\|u_{tt}(t)\right\|^2+\left\|Au_t(t)\right\|^2+\left\|v_t(t)\right\|^2
+\left\|\overline{\eta}^{t}_{t}\right\|^{2}_{L^{2}_{\mu_1}(\mathbb R_+;F_1)}+
\left\|\eta^{t}_{t}\right\|^{2}_{L^{2}_{\mu_2}(\mathbb R_+;F_{1/2})}+\nonumber\\
&+\left\|A^{3/2}u(t)\right\|^2+\omega\left\|A^2u(t)\right\|^2+\omega\left\|Av(t)\right\|^2+\nonumber\\&
+\left\|A^{1/2}v(t)\right\|^{2}+\left\|\overline{T}\overline{\eta}^t\right\|^{2}_{L^{2}_{\mu_1}(\mathbb R_+;F_{1})}+
\left\|T\eta^t\right\|^{2}_{L^{2}_{\mu_2}(\mathbb R_+;F_{1/2})}\leq C^{2}_{R}.\label{mainsmooth}
\end{eqnarray}
   
\begin{lemma} Next statements hold true  
\begin{list}{}{}
\item[(i)] The global attractor $\mathcal{A}$ which existence were established in Section 3 is contained in $D(\mathcal{L})$, the domain of infinitesimal operator $\mathcal{L}$.
\item[(ii)] There exists a positive constant $C_R$ such that for any trajectory $U(t)=(u(t);u_t(t);v(t);\overline{\eta}^t;\eta^t)$ lying in the attractor we have 
\begin{equation}\label{attrest}
\left\|U_t(t)\right\|_{\mathcal{H}}+\left\|\mathcal{L}U(t)\right\|_{\mathcal{H}}\leq C_R,\;\forall t \in \mathbb R.
\end{equation}
\end{list}
\end{lemma}
\rm
{\it Proof.} 

STEP I. Here we use the same ideas as in \cite{ChuiBucci,ChuiLasie,Manuscript}.

Let $\left\{U(t)\equiv(u(t);u_t(t);v(t);\overline{\eta}^t;\eta^t)\right\} \subset \mathcal{H}$ be a full trajectory from the attractor $\mathcal{A}$. Let $\left|\sigma\right|<1$. Applying Theorem \ref{342} with $U^1=U(s+\sigma)$, $U^2=U(s)$ (and, accordingly, the interval $[s,t]$ in place of $[0,t]$), we have that 
\begin{eqnarray} \label{89}
\left\|U(t+\sigma)-U(t)\right\|^{2}_{\mathcal{H}}\leq C_1e^{-\gamma(t-s)}
\left\|U(s+\sigma)-U(s)\right\|^{2}_{\mathcal{H}}+\nonumber\\+
C_2\max\limits_{\tau\in[s,t]}\left\|u(\tau+\sigma)-u(\tau)\right\|^2 
\end{eqnarray}    
for any $t,s \in \mathbb R$ such that $s\leq t$ and for any $\sigma$ with $\left|\sigma\right|<1$. Letting $s\rightarrow -\infty$, \eqref{89} gives
$$
\left\|U(t+\sigma)-U(t)\right\|^{2}_{\mathcal{H}}\leq
C_2\max\limits_{\tau\in(-\infty,t]}\left\|u(\tau+\sigma)-u(\tau)\right\|^2
$$   
for any $t\in\mathbb R$ and $\left|\sigma\right|<1.$ On the attractor we obviously have that 
$$
\frac{1}{\sigma}\left\|u(\tau+\sigma)-u(\tau)\right\|\leq
\frac{1}{\sigma}\int_{0}^{\sigma}\left\|u_t(\tau+t)\right\|dt,\;\;\tau\in\mathbb R.
$$ 

Therefore, by \eqref {87} we obtain that 
$$
\max\limits_{\tau\in \mathbb R}\left\| \frac{U(\tau+\sigma)-U(\tau)}{\sigma}\right\|_{\mathcal{H}} \leq C_R \;\;\mathrm{for}\;\;\left|\sigma\right|<1. 
$$
 
Last estimate implies that function $U(t)$ is absolutely continuous and thus possesses derivative almost everywhere which as well is bounded as follows  
$$
\left\|U_t(t)\right\|_{\mathcal{H}}\leq C_R.
$$

STEP II. Now we prove that $\mathcal{A}\subset D(\mathcal{L})$.
For this we assume that $U_0$ - is a point in the attractor $\mathcal{A}$ that belongs to corresponding full trajectory $\left\{U(t)|t\in \mathbb R\right\}$ that also lies in $\mathcal{A}$ and $U(t)$ possesses a derivative in $t=0$.

Since $U(t)$ is a mild solution of \eqref{7}, then 
$$
U(\sigma)-U_0=e^{\sigma\mathcal{L}}U_0-U_0+\int\limits_{0}^{\sigma}e^{(\sigma-\tau)\mathcal{L}}f(U(\tau))d\tau,\;\;\forall\sigma>0. 
$$          

To check that $U_0$ belongs to the domain of infinitesimal operator $\mathcal{L}$ we need to assure that the following term has a limit as $\sigma\rightarrow 0$ 
$$
\frac{e^{\sigma\mathcal{L}}-I}{\sigma}U_0.
$$

For this we write 
$$
\frac{e^{\sigma\mathcal{L}}-I}{\sigma}U_0=
\underline 
{
\frac{U(\sigma)-U_0}{\sigma}}
-
\underline{\underline{
\frac{1}{\sigma}\int\limits_{0}^{\sigma}e^{(\sigma-\tau)\mathcal{L}}f(U(\tau))d\tau}}.  
$$

Once underlined term converges in force of assumption made in the beginning of step II. We analyse twice underlined term making the following estimate
$$\begin{array}{rl}
\left\|\frac{1}{\sigma}\int\limits_{0}^{\sigma}e^{(\sigma-\tau)\mathcal{L}}\left(f(U(\tau))-f(U_0)\right)d\tau\right\|_{\mathcal{H}}\leq&\frac{1}{\sigma}\int\limits_{0}^{\sigma}\left\|f(U(\tau))-f(U_0)\right\|_{\mathcal{H}}d\tau\leq\\
&\leq\frac{L_R}{\sigma}\int\limits_{0}^{\sigma}\left\|U(\tau)-U_0\right\|_{\mathcal{H}}d\tau\leq\\&\leq
L_R\int\limits_{0}^{\sigma}\left\|\frac{U(\tau)-U_0}{\tau}\right\|_{\mathcal{H}}d\tau\leq L_RC_R\sigma\rightarrow 0.    
\end{array}
$$  

Finally, in view that (see \cite[Theorem 1.2.4.(a)]{Pazy})
$$
\frac{1}{\sigma}\int\limits_{0}^{\sigma}e^{(\sigma-\tau)\mathcal{L}}f(U_0)d\tau\rightarrow f(U_0) \;\;
\mathrm{as}\;\sigma\rightarrow 0,  
$$ 
we make conclusion that $U_0$ belongs to $D(\mathcal{L})$. Using the assertion in Theorem \ref{2.2.1} devoted to classical solutions one can extend the conclusion on whole attractor, thus $\mathcal{A}\subset D(\mathcal{L})$. Besides, it means that the attractor $\mathcal{A}$ consists of full trajectories which correspond to classical solutions of the problem \eqref{7} and then satisfy \eqref{7} literally. It completes the proof of estimate \eqref{attrest}, namely, it gives 
$$
\left\|\mathcal{L}U\right\|_{\mathcal{H}}\leq C_R\;\;\forall U\in \mathcal{A}.
$$  

The proof of the Lemma is complete. 
        
Next Corollary gives more explicit (but not final) form of \eqref{attrest}. For its formulation we set 
$$\begin{array}{rcl}
\phi(t)&\equiv&\beta u(t) + \int\limits_{0}^{+\infty}\mu_1(s)\overline{\eta}^t(s)ds,\\
\rho(t)&\equiv&\phi(t)-\nu A^{-1}v,\\
\psi(t)&\equiv&\omega v(t) + \int\limits_{0}^{+\infty}\mu_2(s)\eta^t(s)ds.
\end{array}
$$
for any classical solution of \eqref{7} $U(t)=(u(t);u_t(t);v(t);\overline{\eta}^t;\eta^t)$.
 
\begin{corollary} \label{smooth2}
There exists a positive constant $C_R$ such that for any trajectory $U(t)=(u(t);u_t(t);v(t);\overline{\eta}^t;\eta^t)$ lying in the attractor we have
\begin{eqnarray}
&\left\|u_{tt}(t)\right\|^2+\left\|Au_t(t)\right\|^2+\left\|v_t(t)\right\|^2
+\left\|\overline{\eta}^{t}_{t}\right\|^{2}_{L^{2}_{\mu_1}(\mathbb R_+;F_1)}+
\left\|\eta^{t}_{t}\right\|^{2}_{L^{2}_{\mu_2}(\mathbb R_+;F_{1/2})}+\nonumber\\
&+\left\|A^{3/2}\phi(t)\right\|^2+\left\|A^2\rho(t)\right\|^2+\left\|A\psi(t)\right\|^2+\nonumber\\
&+\left\|A^{1/2}v(t)\right\|^{2}+\left\|\overline{T}\overline{\eta}^t\right\|^{2}_{L^{2}_{\mu_1}(\mathbb R_+;F_{1})}+
\left\|T\eta^t\right\|^{2}_{L^{2}_{\mu_2}(\mathbb R_+;F_{1/2})}\leq C^{2}_{R}. \label{eplicit}
\end{eqnarray}
\end{corollary}

{\it Proof.} 
First line of \eqref{eplicit} is equivalent to $\left\|U_t\right\|^{2}_{\mathcal{H}}\leq C^{2}_{R}$. Bounds for functions $\rho(t)$ and $\psi(t)$ are obtained directly from original equations \eqref{6}. Remainder of \eqref{eplicit} follows from only 
\begin{equation}\label{forv}
\left\| A^{1/2}v(t)\right\|\leq C_R
\end{equation} if one uses each equality from \eqref{6} again.
In order to prove \eqref{forv} one can repeat procedure of the proof of maximality of the operator $I-\mathcal{L}$ (see Section 2) keeping in mind the goal to obtain \eqref{forv}. 
We propose this method with some insignificant modifications in order to avoid treating with Lax-Millgram Theorem.

It follows from second and fourth equality in \eqref{6} and estimate \eqref{attrest} that 
$$\begin{array}{l}
\omega Av+ \int\limits_{0}^{+\infty}\mu_2(s)A\eta(s)ds=v^*, \;\;\left\|v^*\right\|\leq C_R\\
\eta_s-v=\eta^*,\;\;\left\|\eta^*\right\|_{L^2_{\mu_2}(R_+;F_{1/2})}\leq C_R. 
\end{array}
$$ 

We may integrate second equality and accounting for $\eta(0)=0$ (since $\eta\in D(\mathcal{L})$) we have 
$$
\eta(s)=sv+\int\limits_{0}^{s}\eta^*(y)dy.
$$ 

Now we substitute this to the first equality 
$$
\left\{ \omega+\int\limits_{0}^{+\infty}s\mu_2(s)ds \right\}\cdot Av =
-\int\limits_{0}^{+\infty}\mu_2(s)\int\limits_{0}^{s}A\eta^*(y)dyds+v^*,
$$                   
where right-hand side is obviously estimated by generic constant $C_R$ in space $F_{-1/2}$. 

Thus the proof is complete. 

Using ideas like in Proposition \ref{propVolterra} we are able to continue analysis of attractor's smoothing property. 
 
\begin {lemma} 
There exists a positive constant $C_R$ such that for any trajectory $U(t)=(u(t);u_t(t);v(t);\overline{\eta}^t;\eta^t)$ lying in the attractor we have
\begin{equation} \label{finalsmooth}
\left\|A^{3/2}u(t)\right\|^2+\omega\left\|A^2u(t)\right\|^2+\omega\left\|Av(t)\right\|^2\leq C^{2}_{R}.
\end{equation}
for all $t\in \mathbb R$.
\end{lemma}

{\it Proof.} 
Like in Proposition \ref{propVolterra} we deal with Volterra equation 
\begin{equation}\label{VolterraSmooth}
u(t)-\int\limits_{-\infty}^{t}\frac{\mu_1(t-y)}{\kappa_1+\beta}u(y)dy=h_1(t),
\end{equation}
and due to fully invariance property of $\mathcal{A}$ equality \eqref{VolterraSmooth} holds for all $t\in\mathbb R$ and $h_1(t)\in C(\mathbb R; F_{3/2})$. The same iteration procedure gives $u(t)\in C(\mathbb R;F_{3/2})$ and in addition 
$$
\left\|A^{3/2}u(t)\right\|^2\leq C^{2}_{R}.
$$   

If $\omega>0$ then at first we have to solve 
$$
\omega v+(1-\omega)\int\limits_{-\infty}^{t}k_2(t-y)v(y)dy=h_2(t),\;\;\forall t\in\mathbb R
$$
with $h_2(t)\in C(\mathbb R; F_1)$ and then back to \eqref {VolterraSmooth} with values in $F_2$ instead of $F_{3/2}$.

The proof of the Lemma and inequality \eqref{mainsmooth} is complete. 

\subsection {Exponential attractors.}
In this Subsection we consider sets given by the next definition (according to \cite{Manuscript,Exp})

\begin{definition}
A compact set $A_{exp}\subset\mathcal{H}$ is said to be a {\rm fractal exponential attractor} for the dynamical system $(\mathcal{H},S_t)$ iff $A_{exp}$ is a positively invariant set of finite fractal dimension and for every bounded set $D\subset \mathcal{H}$ there exist positive constants $t_D$, $C_D$ and $\gamma_D$ such that
$$
\sup\limits_{x\in D}\mathrm{dist}_{\mathcal{H}}(S_tx,A_{exp})\leq C_D\cdot e^{-\gamma_D(t-t_D)},\;\;t\geq t_D.
$$ 
\end{definition}

Besides the requirement to be finite dimensional the difference between definition of a global attractor and an exponential attractor is in replacing strict invariance by just positive invariance and in more definite condition on the speed of convergence. The main motivation to consider exponential attractors is that in general case the speed of convergence to the global attractor cannot be estimated. This speed can appear to be small. From the other hand, the exponentiality of the speed to the exponential attractor is guaranteed by the definition.

For the formulation of the Theorem below we introduce an extension of phase space~$\mathcal{H}$ for $\delta>0$
$$
\mathcal{H}_{-\delta}\equiv F_{1-\delta}\times\ F_{-\delta}\times F_{-\delta}\times L^{2}_{\mu_1}(\mathbb R_+;F_{1-\delta})\times L^{2}_{\mu_2}(\mathbb R_+; F_{(1-\delta)/2}).
$$  
  
\begin{theorem} Dynamical system $(\mathcal{H}, S_t)$ possesses a fractal exponential attractor whose dimension is finite in the space~$\mathcal{H}_{-\delta}$, $\delta>0$. 
\end{theorem}
\rm

Proof of the Theorem is based on \cite[Corollary 2.23]{Manuscript} and arguments similar to given in the proof of \cite[Theorem 4.43]{Manuscript}. To provide such arguments we just need to verify that for every $U_0\in\mathcal{B}$ there exists $C_{\mathcal{B},T}$ such that 
\begin{equation}\label{4.131}
\left\|S_{t_1}U_0-S_{t_2}U_0\right\|_{\mathcal{H}_{-\delta}}\leq C_{\mathcal{B},T}\left|t_1-t_2\right|^{\min\left\{\delta,1\right\}},\;\;t_1,t_2\in[0,T],\;U_0\in\mathcal{B},
\end{equation}
where $T>0$ and $\mathcal{B}$ is a positively invariant absorbing set which existence follows from existence of a global attractor and properties of strict Lyapunov function (we may take $\mathcal{B}=\left\{U\in\mathcal{H}|\Phi(U)\leq R\right\}$ for $R>0$ large enough).

Consider $U(t)$ - a classical solution of \eqref{7} with $U_0\in\mathcal{B}$. Then we may estimate (with the help of \eqref{21})
$$
\left\|U_t(t)\right\|_{\mathcal{H}_{-1}}\leq\left\|\mathcal{L}U(t)\right\|_{\mathcal{H}_{-1}}+
\left\|f(U(t))\right\|_{\mathcal{H}_{-1}}\leq C_{\mathcal{B},T}
$$    
and then if $t_1\geq t_2$
$$
\left\|U(t_1)-U(t_2)\right\|_{\mathcal{H}_{-1}}\leq\int\limits_{t_2}^{t_1}\left\|U_t(\tau)\right\|_{\mathcal{H}_{-1}}d\tau \leq C_{\mathcal{B},T}\left|t_1-t_2\right|.  
$$ 

Estimate \eqref{4.131} for $\delta\in(0,1)$ follows from interpolation estimates, e.g., 
$$
\begin{array}{c}
\left\|A^{-\delta}h\right\|\leq\left\|h\right\|^{1-\delta}\left\|h\right\|_{-1}^{\delta},\;h\in H,\\
\int\limits_{0}^{+\infty}\mu_1(s)\left\|A^{1-\delta}\overline{\xi}(s)\right\|ds\leq
\int\limits_{0}^{+\infty}\mu_1(s)^{1-\delta}\left\|A\overline{\xi}(s)\right\|^{1-\delta}\mu_1(s)^{\delta}\left\|\overline{\xi}(s)\right\|^{\delta}ds\leq\\
\leq\left(\int\limits_{0}^{+\infty}\mu_1(s)\left\|A\overline{\xi}(s)\right\|ds\right)^{1-\delta}\left(\int\limits_{0}^{+\infty}\mu_1(s)\left\|\overline{\xi}(s)\right\|ds\right)^{\delta}.
\end{array} 
$$ 

For verification \eqref{4.131} when $U(t)$ is a mild solution we need to approximate $U(t)$ with classical solutions for which \eqref{4.131} has been proved and then pass to limit.

\subsection{Exponential decays to a single equilibrium.}

If the power of the set $\mathcal{N}$ (the set of stationary points) is finite, then conditions of Corollary \ref{finit} hold and each solution of the problem tends to some stationary point (equilibrium point). More actually is true if one imposes some additional conditions, in particular, the speed of convergence to the stationary point might become exponential. 

\begin{definition}
Let an evolution operator $S_t$ be $C^1$ in a Banach space $X$. An equilibrium $e$ is said to be {\rm hyperbolic} if the spectrum $\sigma(L_t)$ of the linear map $L_t=D[S_te]$ satisfies
$$
\sigma(L_t)\cap\left\{z\in\mathbb C: \left|z\right|=1\right\}=\emptyset.
$$ 
for every $t>0$. We also define the index $\mathrm{ind(}e\mathrm{)}$ of the equilibrium $e$ as a dimension of the spectral subspace of the operator $L_1$ corresponding to the set $\sigma_+(L_1)\equiv\left\{z\in\sigma(L_1):\left|z\right|>1\right\}$.
\end{definition}   

Main result of this Subsection relies on the next abstract Theorem (see \cite{ChuiLasie,Manuscript,Raugel} and references therein)

\begin{theorem}  
Let $X$ be a Banach space and the hypotheses of Theorem \ref{31} be in force. Assume that (i) an evolution operator $S_t$ is $C^1$, (ii) the set $\mathcal{N}$ of equilibrium points is finite and all equilibria are hyperbolic, and (iii) there exists a Lyapunov $\mathrm{\Phi}(x)$ function such that 
$$
\mathrm{\Phi}(S_tx)<\mathrm{\Phi}(x), \;\;\forall x\in X,\;\;x\notin\mathcal{N},\; \forall t>0. 
$$   

Then 
\begin{list}{}{}
\item [$\bullet$]For any $y\in X$ there exists $e\in \mathcal{N}$ such that 
$$
\left\|S_ty-e\right\|_{X}\leq C_ye^{-\delta t},\;\;t>0.
$$

Moreover, for any bouded set $B$ in $X$ we have that 
$$
\mathrm{sup}\left\{\mathrm{dist}(S_ty,\mathcal{A})\;:\;y\in B\right\}\leq C_Be^{-\delta t}, \;\;t>0.
$$

Here above $\mathcal{A}$ is a global attractor, $C_y$, $C_B$ and $\delta$ are positive constants, and $\delta$ depends on the minimum, over $e\in \mathcal{N}$, of the distance of the spectrum of $D[S_1e]$ to the unit circle in $\mathbb C$.

\item[$\bullet$]If we assume in addition that (i) $S_1$ is injective on the attractor and (ii)~ the linear map $D[S_1y]$ is injective for every $y\in\mathcal{A}$, then for each $e\in\mathcal{N}$ the unstable manifold $\mathcal{M}^{u}(e)$ is an embedded $C^1$-submanifold of $X$ of finite dimension $\mathrm{ind}(e)$, which implies that $\mathrm{dim}_f\mathcal{A}\leq\max\limits_{e\in\mathcal{N}}\mathrm{ind}(e)$.        
\end{list}\end{theorem} 

Note that all conditions of the Theorem above are verified in corresponding previous subsections except finiteness of the set $\mathcal{N}$ (for discussion of this condition we refer back to Theorem \ref{ChuN} in this article) and hyperbolicity of stationary points. Thus if we consider conditions on stationary points as an assumption we may formulate the following Theorem 

\begin{theorem}
Assume that $\mathcal{N}=\left\{E_i\;:\;i=1,...,n\right\}$ is a finite set. Then the conclusions of Corollary \ref{finit} holds true for the system $(\mathcal{H},S_t)$. In particular, $\mathcal{A}=\cup_{i=1}^{n}\mathcal{M}^{u}(E_i)$. Moreover, if every stationary point is hyperbolic then:
\begin {list}{}{}
\item[$\bullet$] For any $U_0\in\mathcal{H}$ there exists an equilibrium point $E=(e,0,0)\in \mathcal{H}$ and constants $\delta>0$, $C>0$ such that 
$$
\left|S_tU_0-E\right|\leq C_{U_0}e^{-\delta t},\;\; t>0.
$$

Moreover, for any bounded set $B$ in $\mathcal{H}$ we have that 
$$
\sup\left\{\mathrm{dist}\left(S_tU,\mathcal{A}\right)\;:\;U\in B\right\}\leq C_Be^{-\delta t},\;\;t>0.
$$ 

Here above $\mathcal{A}$ is a global attractor, $C_{U_0}$, $C_B$ and $\delta$ are positive constants.  
\item[$\bullet$] For each $E\in \mathcal{N}$ the unstable manifold $\mathcal{M}^{u}(E)$ is an embedded $C^1-$ submanifold of $\mathcal {H}$ of finite dimension $\mathrm{ind}(E)$, which implies that 
$$
\mathrm{dim}_f\mathcal{A}\leq\max\limits_{E\in\mathcal{N}}\mathrm{ind}(E).
$$
\end{list}
\end{theorem}    

\section{Proof of Theorem \ref{342}.} 

The proof of main estimate is based on ideas used in \cite{ChuiLasie} for Von Karman equation. It asserts that a difference of any two solutions can be exponentially stabilized to zero modulo compact perturbation. 

For the sake of reader's convenience we consider the case $\omega=0$ only, which is more complicated. The case $\omega>0$ is simpler because we can use the same representation for nonlinear force as in \cite{ChuiLasie} or \cite{ChuiBucci}.   

Denote 
$$
\kappa_i=\int\limits_{0}^{+\infty}\mu_i(s)ds.
$$

Let $(u^1;v^1;\overline{\eta}^1;\eta^1)$ and $(u^2;v^2;\overline{\eta}^2;\eta^2)$ be two classical solutions of the problem \eqref{7} with initial data $U^i=(u^{i}_{0};u^{i}_{1};v^{i}_{0};\overline{\eta}^{i}_{0};\eta^{i}_{0}),\;i=1,2$ and assume that 
\begin{equation}
\left\|Au^{i}(t)\right\|^2+\left\|u_{t}^{i}(t)\right\|^2+\left\|v^{i}(t)\right\|^2+\left\|\overline{\eta}^{i,t}\right\|_{L^{2}_{\mu_1}(\mathbb R_+;F_1)}^{2}+
\left\|\eta^{i,t}\right\|_{L^{2}_{\mu_2}(\mathbb R_+;F_{1/2})}^{2}\leq R^2
\label{invar}\end{equation}
for $\forall t\geq0$. Also let 
$$
Z(t)\equiv (z(t);z_t(t);\xi(t);\overline{\eta}^t;\eta^t)\equiv
\left(\begin{array}{c}u^{1}(t)-u^{2}(t)\\u_{t}^{1}(t)-u_{t}^{2}(t)\\v^1(t)-v^2(t)\\\overline{\eta}^{1,t}-\overline{\eta}^{2,t}\\\eta^{1,t}-\eta^{2,t}\end{array}\right)^{\mathrm{T}}.
$$

It is clear that components of $Z(t)$ satisfy the equation
\begin{equation}
\left\{
\begin{array}{l}
z_{tt}+\beta A^2z+\int\limits_{0}^{+\infty}\mu_1(s)A^2\overline{\eta}^t(s)ds-\nu A\xi=F(t),\\
\xi_t+\int\limits_{0}^{\infty}\mu_2(s) A\eta^t(s)ds+\nu Az_t=0,\\
\overline{\eta}^{t}_{t}+\overline{\eta}^{t}_{s}=z_t,\;\;\;
\eta_{t}^{t}+\eta_{s}^{t}=\xi,
\end{array}\right.
\label{sysz}\end{equation}     
where 
$$
F(t)=M\left(\left\|A^{1/2}u^2\right\|^2\right)Au^2-M\left(\left\|A^{1/2}u^1\right\|^2\right)Au^1.
$$

To obtain an appropriate form of energy relation from \eqref{sysz} we first transform the term $(F(t),z_t)$. 

\begin{lemma} Let $(u^1(t);v^1(t);\overline{\eta}^{1,t};\eta^{1,t})$ and $(u^2(t);v^2(t);\overline{\eta}^{2,t};\eta^{2,t})$ be classical solutions to problem \eqref{7} satisfying \eqref{invar}. Then following representation 
\begin{equation}
(F(t),z_t)=\frac{d}{dt}Q(t)+P(t)
\label{repr}\end{equation}
holds, where the functions $Q(t)\in C^1(\mathbb R_+)$ and $P(t)\in C(\mathbb R_+)$ satisfy the relations 
\begin{eqnarray}
&&\left|Q(t)\right|\leq C_R \left\|Az\right\|\left\|z\right\| \label{515}\\
&&\left|P(t)\right|\leq C_R \left|\left(\overline{T}\overline{\eta}^{2,t},\overline{\eta}^{2,t}\right)_{L^{2}_{\mu_1}(\mathbb R_+; F_1)}\right|^{1/2}\left(\left\|Az\right\|^2+\left\|z_t\right\|^2\right) \label{516}
\end{eqnarray}  
\end{lemma}
\rm 

{\it Proof.} Introduce the function (the same as in Subsection 2.6)
$$
B(u)=M\left(\left\|A^{1/2}u\right\|^2\right)Au-p.
$$ 

And present $(F(t),z_t(t))$ in following form 
\begin{equation}
(F(t),z_t(t))=\frac{d}{dt}Q_0(t)+P_0(t)
\label{517}\end{equation} 
where 
$$
\begin{array}{l}
Q_0(t)=\int\limits_{0}^{1}\left(B(u^2+\lambda z)-B(u^2),z\right)d\lambda\\
P_0(t)=(B'(u^2)u^{2}_{t},z)-(B(u^1)-B(u^2),u^{2}_{t}) 
\end{array}
$$

Using the differentiability of function $M(z)$ after some straightforward but tedious algebraic manipulations we also have that 
\begin{equation}
P_0(t)=(u_{t}^{2},\mathrm{I}_2\cdot Au^2+\mathrm{I}_1[u^1,u^2]\cdot Az), 
\label{519}\end{equation}
where 
$$
\begin {array}{l} 
\mathrm{I}_{1}[u^1,u^2]=M\left(\left\|A^{1/2}u^1\right\|^2\right)-M\left(\left\|A^{1/2}u^2\right\|^2\right)\\ \\
\mathrm{I}_2=-2M'\left(\left\|A^{1/2}u^2\right\|^2\right)\left(Au^2,z\right)+M\left(\left\|A^{1/2}u^1\right\|^2\right)-M\left(\left\|A^{1/2}u^2\right\|^2\right). 
\end{array}
$$

By using first memory equation we replace $u^{2}_{t}$ appearing in \eqref {519} by 
$$
u^{2}_{t}=\overline{\eta}^{2,t}_{t}+\overline{\eta}^{2,t}_{s}.
$$ 

Substituting this in \eqref{519} written in following form 
$$
P_0(t)=\frac{1}{\kappa_1}\int\limits_{0}^{+\infty}\mu_1(s)(u_{t}^{2},\mathrm{I}_2\cdot Au^2+\mathrm{I}_1[u^1,u^2]\cdot Az)ds  
$$
gives 
$$
P_0(t)=\frac{1}{\kappa_1}\frac{d}{dt}Q_1(t)+\frac{1}{\kappa_1}P_1(t)-\frac{1}{\kappa_1}P_2(t),
$$
with
$$
\begin{array}{l}
Q_1(t)=\int\limits_{0}^{+\infty}\mu_1(s)\left(\overline{\eta}^{2,t}(s),\mathrm{I}_2\cdot Au^2+\mathrm{I}_1[u^1,u^2]\cdot Az\right)ds\\
P_1(t)=\int\limits_{0}^{+\infty}\mu_1(s)\left(\overline{\eta}^{2,t}_{s}(s),\mathrm{I}_2\cdot Au^2+\mathrm{I}_1[u^1,u^2]\cdot Az\right)ds\\
P_2(t)=\int\limits_{0}^{+\infty}\mu_1(s)\left(\overline{\eta}^{2,t}(s),\mathrm{I}_4\cdot Au^2+ \mathrm{I}_2\cdot Au_{t}^{2}+ 2\mathrm{I}_3\cdot Az +\mathrm{I}_1[u^1,u^2]\cdot Az_t\right)ds,
\end{array}
$$
where 
$$
\begin {array}{lcl}
\mathrm{I}_3&=&M'\left(\left\|A^{1/2}u^1\right\|^2\right)\left(Au^1,u^{1}_{t}\right)-M'\left(\left\|A^{1/2}u^2\right\|^2\right)\left(Au^2,u^{2}_{t}\right)  \\
\mathrm{I}_4&=&-4M''\left(\left\|A^{1/2}u^2\right\|^2\right)\left(Au^2,u^{2}_{t}\right)\left(Au^2,z\right)-\\&&-
2M'\left(\left\|A^{1/2}u^2\right\|^2\right)\left[\left(u^{2}_{t},Az\right)+\left(Au^2,z_t\right)\right]+\\&&+
2\left(M'\left(\left\|A^{1/2}u^1\right\|^2\right)\left(Au^1,u^{1}_{t}\right)-M'\left(\left\|A^{1/2}u^2\right\|^2\right)\left(Au^2,u^{2}_{t}\right)\right) 
\end{array}
$$

Thus due to \eqref{517} we have the representation \eqref{repr} with 
\begin{equation}
Q(t)=Q_0(t)+\frac{1}{\kappa_1}Q_1(t)\;\;P(t)=\frac{1}{\kappa_1}(P_1(t)-P_2(t)).
\label{5110}\end{equation}

Now we obtain the estimates for $Q_0(t), Q_1(t), P_1(t)$, and $P_2(t)$. First, let us turn to the analysis of the terms $\mathrm{I}_i$:

1) One can see that: $\left|\mathrm{I}_1[u^1,u^2]\right|\leq C_R\left\|u^1-u^2\right\|$. 

Next representations for terms $\mathrm{I}_i$ allow us to obtain desired estimates:

2) It is straightforward to see that 
$$\begin {array} {l}
\mathrm{I}_2=
\int\limits_{0}^{1}\left[M'\left(\left\|A^{1/2}(u^1-\theta_{\lambda}z)\right\|^2\right)-M'\left(\left\|A^{1/2}u^2\right\|^2\right)\right]d\lambda\cdot\left(Au^2,z\right)+\\
+\int\limits_{0}^{1}\left(\left(M'\left(\left\|A^{1/2}(u^1-\theta_{\lambda}z)\right\|^2\right)-M'\left(\left\|A^{1/2}u^2\right\|^2\right)\right)A(u^1-\theta_{\lambda}z),z\right)d\lambda +\\
+\int\limits_{0}^{1}\left(M\left(\left\|A^{1/2}u^2\right\|^2\right)\left(A(u^1-\theta_{\lambda}z)-Au^2\right),z\right)d\lambda +\\+
\int\limits_{0}^{1}M'\left(\left\|A^{1/2}(u^1-\theta_{\lambda}z)\right\|^2\right)\left(\theta_{\lambda}Az,z\right)d\lambda  
\end {array}
$$ 
where $\theta_{\lambda}\in(0,1)$ satisfies the equality: 
$$
\left\|A^{1/2}(u^1-\theta_{\lambda}z)\right\|^2=(1-\lambda)\left\|A^{1/2}u^2\right\|^2+\lambda\left\|A^{1/2}u^1\right\|^2
$$

Hence, $\left|\mathrm{I}_2\right|\leq C_R\left\|Az\right\|\left\|z\right\|$.

3)It is elementary to see that 
$$\begin {array} {lrr}
\mathrm{I}_3&=&\left[M'\left(\left\|A^{1/2}u^1\right\|^2\right)-M'\left(\left\|A^{1/2}u^2\right\|^2\right)\right]\left(Au^1,u^{1}_{t}\right)+\\&&+
M'\left(\left\|A^{1/2}u^2\right\|^2\right)\left[\left(Az,u_{t}^{1}\right)+\left(Au^2,z_t\right)\right]
\end{array} 
$$

Hence, $\left|\mathrm{I}_3\right|\leq C_R\left(\left\|Az\right\|+\left\|z_t\right\|\right)$.
  
4)One can also see that
$$\begin {array}{rll}
\mathrm{I}_4&=&2\left(Au^2,u^{2}_{t}\right)\mathrm{I}_{2}^{*}+2M'\left(\left\|A^{1/2}u^2\right\|^2\right)\left(Az,z_t\right)+\\&&+
\left(M'\left(\left\|A^{1/2}u^1\right\|^2\right)-M'\left(\left\|A^{1/2}u^2\right\|^2\right)\right)
\left(\left(Az,u^{1}_{t}\right)+\left(Au^2,z_t\right)\right)
\end {array}
$$
where 
$$
\mathrm{I}_{2}^{*}=-2M''\left(\left\|A^{1/2}u^2\right\|^2\right)\left(Au^2,z\right)+M'\left(\left\|A^{1/2}u^1\right\|^2\right)-M'\left(\left\|A^{1/2}u^2\right\|^2\right)
$$

Note that $\mathrm{I}_{2}^{*}$ admits the same estimate as $\mathrm{I}_2$.

Hence, $\left|\mathrm{I}_4\right|\leq C_R\left(\left\|Az\right\|^2+\left\|z_t\right\|^2\right)$. 

Now we are able to prove necessary bounds pertaining to the terms $Q_0(t), Q_1(t)$, $P_1(t)$, and $P_2(t)$. Since  
$$\begin {array} {l}
Q_0(t)=-\int\limits_{0}^{1}\mathrm{I}_1[u^2+\lambda z,u^2]d\lambda\left(Au^2,z\right)-
\int\limits_{0}^{1}\lambda M\left(\left\|A^{1/2}(u^2+\lambda z)\right\|^2\right)d\lambda\left(Az,z\right)   
\end{array}
$$
we obviously have that $\left|Q_0(t)\right|\leq C_R\left\|Az\right\|\left\|z\right\|$.

Using the expressions of $Q_1(t)$, $P_1(t)$, $P_2(t)$ and estimates for $\mathrm{I}_i$ we obtain other inequalities:  
$$\begin{array}{l}
\left|Q_1(t)\right|\leq C_R\left\|Az\right\|\left\|z\right\| \\
\left|P_1(t)\right|\leq C_R\left|\left(\overline{T}\overline{\eta}^{2,t},\overline{\eta}^{2,t}\right)_{L^{2}_{\mu_1}(\mathbb R_+; F_1)}\right|^{1/2}\left\|Az\right\|^2\\
\left|P_2(t)\right|\leq C_R\left\|\overline{\eta}^{2,t}\right\|_{L^{2}_{\mu_1}(\mathbb R_+; F_1)}\left(\left\|Az\right\|^2+\left\|z_t\right\|^2\right).
\end{array}$$

Estimate for $P_1(t)$ were obtained in view of the following observation. Consider any $w\in H$, then
$$\begin{array}{c}
\left|\int_{0}^{+\infty}\mu_1(s)(\overline{\eta}^{2,t}_{s},w)ds\right|\leq
\int\limits_{0}^{+\infty}(-\mu'_1(s))\left\|\overline{\eta}^{2,t}\right\|_1ds\cdot\left\|A^{-1}w\right\|\leq\\ \leq\mu_{1}^{1/2}(0)\left(\int\limits_{0}^{+\infty}(-\mu'_1(s))\left\|\overline{\eta}^{2,t}\right\|_{1}^{2}ds\right)^{1/2}\left\|A^{-1}w\right\|\leq\\\leq \mu_{1}^{1/2}(0)\left|\left(\overline{T}\overline{\eta}^{2,t},\overline{\eta}^{2,t}\right)_{L^{2}_{\mu_1}(\mathbb R_+; F_1)}\right|^{1/2}\left\|A^{-1}w\right\|.
\end{array}$$

The final estimate is derived in view that $$
\left\|\overline{\eta}^{2,t}\right\|_{L^{2}_{\mu_1}(\mathbb R_+; F_1)}\leq \left|\left(\overline{T}\overline{\eta}^{2,t},\overline{\eta}^{2,t}\right)_{L^{2}_{\mu_1}(\mathbb R_+; F_1)}\right|^{1/2}.$$ 

The proof of Lemma is complete.

{\bf Proper proof of Theorem \ref{342}}

By \eqref{repr} for these solutions we have energy relation 
\begin{equation}
\frac{d}{dt}\mathcal{E}^{0}(t)=\left(\overline{T}\overline{\eta}^t,\overline{\eta}^t\right)_{L^{2}_{\mu_1}(\mathbb R_+;F_1)}+\left(T\eta^t,\eta^t\right)_{L^{2}_{\mu_2}(\mathbb R_+;F_{1/2})}+P(t)
\label{5111}\end{equation}
where 
$$\begin{array}{c}
\mathcal{E}^{0}(t)=\frac{1}{2}\left[\left\|z_t(t)\right\|^2+\left\|Az(t)\right\|^2+\left\|\xi(t)\right\|^2\right]+\\+\frac{1}{2}\left[\left\|\overline{\eta}^t\right\|^{2}_{L^{2}_{\mu_1}(\mathbb R_+;F_1)}+\left\|\eta^t\right\|^{2}_{L^{2}_{\mu_2}(\mathbb R_+;F_{1/2})}-2Q(t)\right].
\end{array}$$

It follows from \eqref{515} that 
\begin{equation}
\frac{3}{8}\left\|Z(t)\right\|_{\mathcal{H}}^{2}-C_R\left\|z(t)\right\|^2\leq
\mathcal{E}^{0}(t)\leq
\frac{5}{8}\left\|Z(t)\right\|_{\mathcal{H}}^{2}+C_R\left\|z(t)\right\|^2
\label{5112}\end{equation}

Now we consider 
\begin{equation}
V(t)\equiv \mathcal{E}^{0}(t)+\sum\limits_{i=1}^{3}\varepsilon_i\Phi_{i}(t)
\label{5113}\end{equation}
where 
$$\begin{array}{l}
\Phi_1(t)=(z_t,z)\\
\Phi_2(t)=-(A^{-2}z_t,\overline{\eta}^t)_{L^{2}_{\mu_1}(\mathbb R_+;F_{1})}\\
\Phi_3(t)=-(\nu z+A^{-1}\xi,\eta^t)_{L^{2}_{\mu_2}(\mathbb R_+;F_{1/2})}.
\end{array}
$$

Positive constants $\varepsilon_i$ will be chosen in the sequel. For $V(t)$ we have estimate similar to \eqref{5112}
\begin{equation}
\frac{1}{4}\left\|Z(t)\right\|_{\mathcal{H}}^{2}-C_R\left\|z(t)\right\|^2\leq
V(t)\leq
\left\|Z(t)\right\|_{\mathcal{H}}^{2}+C_R\left\|z(t)\right\|^2
\label{5114}\end{equation} 
as soon as common sum of $\varepsilon_i$ is sufficiently small. 

Now we compute derivatives of $\Phi_i(t)$
$$
\begin{array}{rcl}
\Phi'_{1}(t)&=& (z_{tt},z)+\left\|z_t\right\|^2=\\&=&\left(-\beta A^2z-\int\limits_{0}^{+\infty}\mu_1(s)A^2\overline{\eta}^t(s)ds+\nu A\xi+F(t),z\right)+\left\|z_t\right\|^2=\\
&=&-\beta\left\|Az\right\|^2-\int\limits_{0}^{+\infty}\mu_1(s)(\overline{\eta}^t(s),z)_1ds+\nu(\xi,Az)+(F(t),z)+\left\|z_t\right\|^2\\
\Phi'_{2}(t)&=& -(A^{-2}z_{tt},\overline{\eta}^t)_{L^{2}_{\mu_1}(\mathbb R_+;F_{1})}-
(A^{-2}z_t,-\overline{\eta}^{t}_{s}+z_t)_{L^{2}_{\mu_1}(\mathbb R_+;F_{1})}=\\
&=&\int\limits_{0}^{+\infty}\mu_1(s)\left(\beta A^2z +\int\limits_{0}^{+\infty}\mu_1(\tau)A^2\overline{\eta}^t(\tau)d\tau-\nu A\xi-F(t),\overline{\eta}^{t}(s)\right)ds+\\&&+\int\limits_{0}^{+\infty}\mu_1(s)(z_t,\overline{\eta}^{t}_{s})ds-\kappa_1\left\|z_t\right\|^2\\
\Phi'_{3}(t)&=&-(\nu z_t+ A^{-1}\xi_t,\eta^t)_{L^{2}_{\mu_2}(\mathbb R_+;F_{1/2})}-(\nu z+A^{-1}\xi,\eta^t)_{L^{2}_{\mu_2}(\mathbb R_+;F_{1/2})}=\\
&=&(\int\limits_{0}^{+\infty}\mu_2(\tau)\eta^{t}(\tau)d\tau,\eta^t)_{L^{2}_{\mu_2}(\mathbb R_+;F_{1/2})}-\nu\kappa_2(Az,\xi)-\kappa_2\left\|\xi\right\|^2+\\&&+(\nu z+A^{-1}\xi,\eta^{t}_{s})_{L^{2}_{\mu_2}(\mathbb R_+;F_{1/2})}.
\end{array}
$$

Our main task is to estimate the term $\frac{d}{dt}V(t)+\gamma \left\|Z(t)\right\|^{2}_{\mathcal{H}}$ with small parameter $\gamma>0$, that could be chosen in next steps of the proof, by the the sum of next form 
$$
-\alpha\left\|Z(t)\right\|^{2}_{\mathcal{H}}+P(t)+C_R\left\|z(t)\right\|^2.
$$ 

For this we rewrite inequality for $\frac{d}{dt}\mathcal{E}^{0}(t)$, via
$$\begin{array}{c}
\frac{d}{dt}\mathcal{E}^0(t)\leq -\frac{\delta}{4}\left\|\overline{\eta}^t\right\|^{2}_{L^{2}_{\mu_1}(\mathbb R_+;F_{1})}- \frac{\delta}{4}\left\|{\eta}^t\right\|^{2}_{L^{2}_{\mu_2}(\mathbb R_+;F_{1/2})}-\\-
\frac{1}{4}\int\limits_{0}^{+\infty}(-\mu'_1(s))\left\|\overline{\eta}^t(s)\right\|^{2}_{1}ds-
\frac{1}{4}\int\limits_{0}^{+\infty}(-\mu'_2(s))\left\|{\eta}^t(s)\right\|^{2}_{1/2}ds+P(t).
\end{array}
$$

Further steps contain splitings of scalar products according to Coushy inequality. We may choose $\varepsilon_i$ small enough for all products where memory variables are included to be splitted in such way that terms of the form $\left\|Az\right\|^2,\left\|z_t\right\|^2,\left\|\xi\right\|^2$ won't give an essential contribution to the general estimate, for example 
$$
\varepsilon_1\int\limits_{0}^{+\infty}\mu_1(s)(z,\overline{\eta}^t(s))_1ds\leq
\frac{\varepsilon_1}{2\sigma}\left\|\overline{\eta}\right\|^{2}_{L^{2}_{\mu_1}(\mathbb R_+;F_{1})}+\varepsilon_1\frac{\sigma}{2} \left\| Az\right\|^2,\;\;\forall\sigma>0.  
$$       

Here we first need to choose small enough $\sigma$ (for the coefficient near $\left\| Az\right\|^2$) and then $\varepsilon_1$ (for the one near $\left\|\overline{\eta}\right\|^{2}_{L^{2}_{\mu_1}(\mathbb R_+;F_{1})}$). Because of the presence of terms with derivatives with the respect to $s$ (for instance, $\overline{\eta}^{t}_{s}$) we picked out terms with $\mu'_i$ in the inequality for $\frac{d}{dt}\mathcal{E}^{0}(t)$. Now we vanish the coefficient near $(Ax,\xi)$, for this we set $\varepsilon_1=\kappa_2\varepsilon_3$. Besides, the setting $\varepsilon_2=\frac{2}{\kappa_1}\varepsilon_1$ gives negative coefficient near $\left\|z_t\right\|^2$. 

Finally, $(F(t),z(t))\leq \frac{\sigma}{2}\left\|Az(t)\right\|^2+\frac{1}{2\sigma}\left\|z(t)\right\|^2$ for all $\sigma>0$. 
Furthermore, due to \eqref{5113} we may choose small enough $\gamma>0$ such as 
$$
\frac{d}{dt}V(t)+\gamma V(t)\leq C_R\left\|z(t)\right\|^2+C_R\left|\left(\overline{T}\overline{\eta}^{2,t},\overline{\eta}^{2,t}\right)_{L^{2}_{\mu_1}(\mathbb R_+; F_1)}\right|\left(\left\|Az\right\|^2+\left\|z_t\right\|^2\right) 
$$     

Here we again used Coushy inequality to obtain 
$$
\left|\left(\overline{T}\overline{\eta}^{2,t},\overline{\eta}^{2,t}\right)_{L^{2}_{\mu_1}(\mathbb R_+; F_1)}\right|\left(\left\|Az\right\|^2+\left\|z_t\right\|^2\right)
$$ 
instead of 
$$
\left|\left(\overline{T}\overline{\eta}^{2,t},\overline{\eta}^{2,t}\right)_{L^{2}_{\mu_1}(\mathbb R_+; F_1)}\right|^{1/2}\left(\left\|Az\right\|^2+\left\|z_t\right\|^2\right).
$$

After using Gronwall Lemma we obtain 
$$
\begin{array}{c}
\left\|Z(t)\right\|^{2}_{\mathcal{H}}\leq C_R\left\|Z(0)\right\|^{2}_{\mathcal{H}}e^{-\gamma t}+C_R\max\limits_{\tau\in[0,t]}\left\|z(t)\right\|^2+\\
+C_R\int\limits_{0}^{t}e^{-\gamma(t-\tau)}\left|\left(\overline{T}\overline{\eta}^{2,\tau},\overline{\eta}^{2,\tau}\right)_{L^{2}_{\mu_1}(\mathbb R_+; F_1)}\right|\left\|Z(\tau)\right\|^{2}_\mathcal{H}d\tau
\end{array}
$$

Now using the fact 
$$
\int\limits_{0}^{+\infty}\left|\left(\overline{T}\overline{\eta}^{2,t},\overline{\eta}^{2,t}\right)_{L^{2}_{\mu_1}(\mathbb R_+; F_1)}\right|dt\leq C_R,
$$
which follows from energy relation and inequality \eqref{invar}, and Gronwall Lemma of the form of Lemma \ref{A2} (see below) setting 
$$
\begin {array}{c}
\phi(t)=\left\|Z(t)\right\|^{2}_{\mathcal{H}}e^{\gamma t},\;\;\;\;\phi_1(t)=\left\|Z(0)\right\|^{2}_{\mathcal{H}}+C_Re^{\gamma t}\max\limits_{\tau \in[0,t]}\left\|z(\tau)\right\|^2,\\ 
\phi_2(t)=\left|\left(\overline{T}\overline{\eta}^{2,t},\overline{\eta}^{2,t}\right)_{L^{2}_{\mu_1}(\mathbb R_+; F_1)}\right|.
\end{array}
$$ 
we obtain stabilizability estimate.

\begin{lemma} \label{A2}Let $\phi(t),\;\phi_1(t)$ and $\phi_2(t)$ be scalar positive functions. We also assume that $\phi_1$ is a non-decreasing function and $\phi_2$ satisfies the following condition
\begin{equation}
\int\limits_{0}^{+\infty}\phi_2(t)dt<\infty.
\label{david}\nonumber \end{equation}

Besides, the relation 
\begin{equation}
\phi(t)\leq\phi_1(t)+C_1\int\limits_{0}^{t}\phi_2(\tau)\phi(\tau)d\tau
\label{smile}\nonumber \end{equation}
holds for all $t\geq0$. Then there exists positive constant $C$ such as 
$$
\phi(t)\leq C\phi_1(t)\;\; \forall t\geq 0.
$$\end{lemma}
\rm 

\end{document}